%% file: DGFEM-Linear-Dirichlet-Curved.tex
\documentclass[twoside,a4paper]{article} 
\usepackage{a4wide} 
\usepackage{hyperref}
\usepackage{latexsym} 
\usepackage{amssymb} 
\usepackage{amsmath} 
\usepackage{theorem} 
\usepackage{mathrsfs} 
\usepackage{xspace}  
\usepackage{fancyhdr} 
\usepackage{float} 
\usepackage{listings}
\usepackage{enumerate}
\usepackage{pdfpages}
\usepackage{morefloats}
\usepackage{tikz}
\usepackage{pdfpages}
\usepackage{accents}
\usepackage{graphicx} 
\usepackage{caption}
\usepackage{rotating}
\usepackage{subcaption}
\definecolor{SussexFlint}{rgb}{.00,.19,.21}
\definecolor{SussexGrey}{rgb}{.51,.58,.49}
\definecolor{SussexOrange}{rgb}{.94,.29,.00}
\definecolor{SussexYellow}{rgb}{1.00,.73,.00}
\definecolor{SussexRed}{rgb}{.94,.01,.49}
\definecolor{SussexPurple}{rgb}{.48,.06,.44}
\definecolor{SussexGreen}{rgb}{.00,.58,.46}
\definecolor{SussexBlue}{rgb}{.00,.58,.65}
\colorlet{a}{SussexOrange}
\colorlet{b}{SussexRed}
\colorlet{c}{SussexYellow}
\colorlet{d}{SussexPurple}
\colorlet{e}{SussexGreen}
\colorlet{f}{SussexBlue}
\colorlet{g}{SussexGrey}
\colorlet{h}{white}
\colorlet{i}{black}
\colorlet{j}{SussexFlint}
\numberwithin{equation}{section}
\newtheorem{theorem}{Theorem}[section]
\newtheorem{lemma}[theorem]{Lemma}
\newtheorem{corollary}[theorem]{Corollary}
\newtheorem{definition}[theorem]{Definition}%
\newtheorem{remark}[theorem]{Remark}
\newcommand{\qed}{\hfill$\square$}

\newtheorem{example}[theorem]{Example}

\newcommand{\average}[1]{\ensuremath{\langle\!\langle#1\rangle\!\rangle} }
\newcommand{\avg}[1]{\ensuremath{\langle\!\langle#1\rangle\!\rangle} }
\newcommand{\jump}[1]{\ensuremath{[\![#1]\!]} }

\newcommand{\Th}{\ensuremath{\mathscr{T}_h}}
\newcommand{\Ta}{\ensuremath{\mathbf{T}}}
\newcommand{\mc}{\ensuremath{\mathcal{H}}}

\newcommand{\Eb}{\ensuremath{\mathscr{E}^b_h}}
\newcommand{\Ei}{\ensuremath{\mathscr{E}^i_h}}
\newcommand{\Eib}{\ensuremath{\mathscr{E}^{i,b}_h}}

\newcommand{\deter}{\ensuremath{\operatorname{det}}}

\newcommand{\gradTa}{\ensuremath{\nabla_\Ta}}

\newcommand{\Reals}{\ensuremath{\mathbb{R}}}

\title{A DGFEM for Nondivergence Form Elliptic Equations with Cordes Coefficients on Curved Domains}
\author{Ellya Kawecki\footnote{EK acknowledges support of the Engineering and Physical Sciences Research Council
[EP/L015811/1].}}
\begin{document}
\maketitle
\begin{abstract}
In ``I. Smears, E. S\"{u}li, \emph{Discontinuous Galerkin finite element approximation of nondivergence form elliptic equations with Cord\'{e}s coefficients. SIAM J. Numer Anal., 51(4):2088-2106, 2013}"
the authors designed and analysed a discontinuous Galerkin finite element method for the approximation of solutions to elliptic partial differential equations in nondivergence form. The results were proven, based on the assumption that the computational domain was convex and \emph{polytopal}. In this paper, we extend this framework, allowing for Lipschitz continuous domains with piecewise curved boundaries. 
\end{abstract}
\begin{section}{Introduction}
In this paper, we also tackle the problem of approximating solutions to elliptic problems in \emph{nondivergence} form on domains with curved boundaries. That is, the problems that fall into this framework do not in general possess a weak formulation; which means that the standard finite element approach (which is centred around a weak formulation) is no longer justified. Furthermore, unlike elliptic equations in divergence form, uniform ellipticity is not, in general, a sufficient assumption for well--posedness of the corresponding boundary value problem, when the coefficients are not continuous, see~\cite{MR1814364,MR2260015} for such examples.

We extend the framework found in~\cite{MR3077903} to problems with domains that are allowed to have piecewise curved boundaries, under the assumption that the curved portions of the boundary have nonnegative curvature. We note, however, that this does not restrict the framework in any way, since the scheme we define coincides with the scheme introduced in~\cite{MR3077903} when the domain is polytopal. 

One can intuitively characterise the difference between convexity and uniform convexity (indeed the latter implies the former), as follows: a domain is convex if the line segment between any two points in the domain is contained in the domain. A domain is uniformly convex if the former statement is true, and if the open line segment between any two points on the boundary of the domain is contained in the interior of the domain. For example, the unit square $(0,1)^2$ is convex, but not uniformly convex. In contrast, the unit disc $\{(x,y):x^2+y^2<1\}$ is uniformly convex.

As will be seen throughout this paper, extending  the framework of~\cite{MR3077903} is nontrivial, both in terms of reformulating the numerical method and adapting the consistency, stability and error analysis of this new method. In particular, there are new terms arising in the bilinear form (due to the curvature of the boundary), which are necessary to maintain the consistency of the method; we are able to demonstrate this necessity in Section~\ref{Experiments}, by attempting to implement the method present in~\cite{MR3077903}, without our adaptation, to an elliptic problem with a curved boundary. In this experiment, we will see both a lack of consistency, and error results inferior to those produced by the new method we propose.

Elliptic equations in nondivergence form appear in applications to fully nonlinear second order elliptic partial differential equations (PDEs), which, upon employing a suitable linearisation technique (Newton's method for instance), produces an infinite sequence of such equations. Two particular examples of  nonlinear problems are the equations of Monge--Amp\`ere (MA) and Hamilton--Jacobi--Bellman (HJB) type. The MA problem arises in areas such as optimal transport and differential geometry, and has been an area of interest, both from an analytical and a numerical computation point of view for many years,  see~\cite{MR1426885,MR0180763,MR1454261,MR1624426} and~\cite{MR2916369,MR3162358,MR3653852}; while the HJB problem arises in applications to engineering, physics, economics and finance~\cite{MR2179357}, where~\cite{MR3033005,MR3196952} mark recent  developments in the numerical analysis of such problems. 

The paper is organised as follows: In Section~\ref{The PDE} we introduce the type of equations under consideration in more detail, and provide an important existence and uniqueness result for such problems. In Section~\ref{The numerical method} we begin by introducing the notation needed, we define our numerical method, provide the necessary structural assumptions for the approximation of the computational domain, and then define the numerical method. In Section~\ref{Analysis of the numerical method} we prove a stability result for one of the main bilinear terms our numerical scheme; this stability result is then used as a main tool in the proof of existence and uniqueness of a numerical solution, and proceed to prove an important consistency result. Finally we prove an error estimate that is optimal in terms of the mesh size. In Section~\ref{Experiments}, we run several numerical experiments where the true solution is known. This allows us to verify the error estimates presented in Section~\ref{Analysis of the numerical method}, and test the robustness of the scheme by considering operators with discontinuous coefficients, as well as a nonconvex computational domain. Section~\ref{Conclusion} is the final section, where we give concluding remarks on what has been accomplished in this paper, as well as plans for future research.

\subsection{The PDE}\label{The PDE}
Consider the following second-order elliptic boundary--value problem: find $u:\Omega\to\mathbb R$ such that
\begin{equation}\label{1}
\left\{
\begin{aligned}
Lu&= f,\quad\mbox{in}\quad\Omega,\\
u  &= 0,\quad\mbox{on}\quad\partial\Omega,
\end{aligned}
\right.
\end{equation}
where $\Omega\subset\mathbb R^d$ is a Lipschitz continuous domain, and
$$Lu:=\sum_{i,j=1}^d A_{ij}D^2_{ij}u,$$ where $A\in L^\infty(\Omega;\mathbb R^{d\times d})$, is uniformly elliptic, and $f\in L^2(\Omega)$ is a given function. Furthermore, we assume that $A$ satisfies the \emph{Cordes condition}: there is an $\varepsilon\in(0,1]$ such that
\begin{equation}\label{cordes}
\frac{|A|^2}{(\operatorname{Tr}(A))^2}=\frac{\sum_{i,j=1}^dA_{ij}^2}{(\sum_{i=1}^dA_{ii})^2}\le\frac{1}{d-1+\varepsilon}\quad\mbox{a.e. in }\Omega.
\end{equation}
We quote the following result from~\cite{MR3077903}, Theorem 3.
\begin{theorem}\label{Thm:1}
Assume that $\Omega$ is convex, and that $A$ satisfies~(\ref{cordes}). Then, there exists a unique $u\in H^2(\Omega)\cap H^1_0(\Omega)$ that is a strong solution of~(\ref{1}).
\end{theorem}
\end{section}
\begin{section}{The numerical method}\label{The numerical method}
In this section we detail the numerical method used to approximate solutions of a renormalised, equivalent representation of~(\ref{1}). To this end, we consider the renormalisation function  $\gamma:\Omega\to\mathbb R^+$, defined by
$$\gamma:=\frac{\operatorname{Tr}A}{|A|^2}.$$
 Taking into account that $\gamma$ is uniformly positive (when restricted to matrix--valued functions with uniformly positive trace), we see that
 $u\in H^2(\Omega)\cap H^1_0(\Omega)$ solves 
\begin{equation}\label{eqwgamma}
\left\{
\begin{aligned}
\gamma Lu&= \gamma f,\quad\mbox{in}\,\,\Omega,\\
u  &= 0,\quad~~\mbox{on}\,\,\partial\Omega,
\end{aligned}
\right.
\end{equation}
uniquely, if and only if it is the unique solution of~(\ref{1}).
\subsection{Finite element spaces}
The finite element spaces we consider consist of discontinuous piecewise polynomial functions, and as such we must define some important notation. 

Let $\{\Th\}_h$ be a sequence of shape-regular meshes
on $\Omega$, consisting of simplices or parallelepipeds, possibly with curved edges.

\emph{Edge and vertex sets.}
Given a triangulation $\Th$, we denote by $\mathscr{E}_h$, the set of edges of $\Th$, by $\mathscr{E}_h^i$ the set of interior edges of $\Th$, by $\mathscr{E}^b$, the set of boundary edges, by $\mathscr{E}_h^{i,b}:=\mathscr{E}_h^{i}\cup\mathscr{E}_h^{b}$, and by $\mathscr{V}^b_h$ the set of boundary vertices of $\Th$.

\emph{Exact approximation.} We assume that our mesh approximates the domain exactly. That is,
\begin{equation}\label{Exact}
\bigcup_{K\in\mathscr{T}_h}\overline{K}=\overline{\Omega};
\end{equation}
this assumption is considered in the sense of~\cite{MR1014883}. The (possibly curved) open simplices $K$ are the images of a reference open simplex $\hat K$ under a collection of mappings
$$F_K=\tilde{F}_K+\Phi_K,$$
where 
\begin{equation}\label{def:affine-map}
\tilde{F}_K:\hat{x}\mapsto\tilde{B}_K\hat x+\tilde{b}_K
\end{equation}
is an invertible affine map and $\Phi_K\in C^1(\hat{K};\mathbb R^d)$ satisfies
$$C_K:=\sup_{\hat{x}\in\hat K}\|D\Phi_K(\hat x)\tilde{B}_K^{-1}\|<1,$$
where $\|\cdot\|$ denotes the Euclidean norm on $\mathbb R^d$, when the entry is vector valued, and the induced Euclidean norm, when the entry is matrix valued.
\begin{remark}
The analysis of our numerical method does not require the domain, $\Omega$, to be convex. This means that there are cases where a unique solution to our numerical method exists, but the assumptions of Theorem~\ref{Thm:1} are not satisfied. This motivates the following two definitions.
\end{remark}
\begin{definition}[Piecewise $C^{k,\alpha}$ domain]
A domain $\Omega\subset\mathbb R^d$ is piecewise $C^{k,\alpha}$ for $k\in\mathbb N$, $\alpha\in(0,1)$, if we may express the boundary of $\Omega$, $\partial\Omega$, as a finite union
\begin{equation}\label{piecewiserep}
\partial\Omega = \bigcup_{n=1}^N\overline{\Gamma}_n,
\end{equation}
where each $\Gamma_n\subset\mathbb R^d$ is of zero $d$-dimensional Lebesgue measure, and admits a local representation as the graph of a $C^{k,\alpha}$ function.
\end{definition}
\begin{definition}[Piecewise $C^{k,\alpha}$-piecewise convex domain]\label{p-conv-def}
A domain $\Omega\subset\mathbb R^d$, is a piecewise $C^{k,\alpha}$-piecewise convex domain, with $k\ge 2$, $\alpha\in(0,1)$, if $\Omega$ is Lipschitz continuous, and the boundary of $\Omega$, $\partial\Omega$, is given by a finite union of the form~(\ref{piecewiserep}), and each $\Gamma_n\subset\mathbb R^d$, can be expressed as a subset of a level set of a $C^{k,\alpha}$ \emph{convex} function $\varphi_n$. Moreover, the unit normal to $\Gamma_n$, $n_{\Gamma_n}$, must correspond to 
\begin{equation}\label{convex-normal}
n_{\Gamma_n} = \frac{\nabla\varphi_n}{|\nabla\varphi_n|}.
\end{equation}
\end{definition}
\begin{remark}
Note that if $\Omega$ is a piecewise $C^{k,\alpha}$-piecewise convex domain, with $k\ge 2$, $\alpha\in(0,1)$, it is not in necessarily convex. For example, consider the ``key-hole shaped" domain 
\begin{equation}\label{Key-hole}
\Omega = \{x^2+y^2<1:y\ge1/\sqrt{2}\}\cup[-1/\sqrt{2},1/\sqrt{2}]\times[-3,1/\sqrt{2}].
\end{equation}
See Figure~\ref{key:hole:mesh}.
\end{remark}
\begin{remark}
The unit normal assumption~(\ref{convex-normal}), in Definition~\ref{p-conv-def}, is required to exclude piecewise smooth domains with boundary portions of strictly negative curvature, for example the subset of $\mathbb R^2$ given by:
\begin{equation}\label{diskremoved}
([-2,0]\times[0,2])\setminus\{x^2+y^2<1/4\}.
\end{equation}
See Figure~\ref{key:hole:mesh}.
\end{remark}

\emph{Mesh size.} For each element $K\in\Th$, let $h_K:=\operatorname{diam}\tilde{K}\ge C(d)\|\tilde{B}_K\|$ (where $\tilde{K}=\tilde{B}_K(\hat K)$). It is assumed that $h = \max_{K\in\Th} h_K$ for each mesh $\Th$.

\emph{Mesh conditions.} We shall adopt the following assumptions on the meshes. The meshes are allowed to be irregular, i.e., there may be hanging nodes. We assume that there is a uniform upper bound on the number of edges composing the boundary of any given element; in other words, there is a $c_\mathcal{F}>0$, independent of $h$, such that
\begin{equation}\label{meshcond1}
\max_{K\in\mathscr{T}_h}\operatorname{card}\{F\in\mathscr{E}^{i,b}_h:F\subset\partial K\}\le c_\mathcal{F}\quad\forall K\in\mathscr{T}_h,~\forall h>0.
\end{equation}
It is also assumed that any two elements sharing an edge have commensurate diameters, i.e., there is a $c_{\mathcal{T}}\ge 1$, independent of $h$, such that
\begin{equation}\label{meshcond2}
\max(h_K,h_{K'})\le c_\mathcal{T}\min(h_K,h_{K'}),
\end{equation}
for any $K$ and $K'$ in $\mathscr{T}_h$ that share an edge. 
Furthermore, for each edge $F\in\mathscr{E}^{i,b}_h$, we define
\begin{equation}
\tilde{h}_F:=\left\{\begin{array}{l l}
\min(h_K,h_{K'}) & \mbox{if}~F\in\mathscr{E}^i_h,\\
h_K & \mbox{if}~F\in\mathscr{E}^b_h.
\end{array}\right.
\end{equation}
where $K$ and $K'$ are such that $F=\partial K\cap\partial K'$ if $F\in\mathscr{E}^i_h$, or $F\subset\partial K\cap\partial\Omega$ if $F\in\mathscr{E}^b_h$. 
The assumptions on the mesh, 
in particular~(\ref{meshcond2}), show that if $F$ is an edge of $K$, then
\begin{equation}
h_K\le c_\mathcal{T}\tilde{h}_F.
\end{equation}
Finally, we assume that each $F\in\Eb$ satisfies
\begin{equation}\label{Fcontained}
F = F\cap\Gamma_n,
\end{equation}
for some $n\in\{1,\ldots,N\}$. This implies that each boundary face is completely contained in a boundary portion $\Gamma_n$.
\begin{remark}
In order to prove the error estimate in Section~\ref{Error estimates}, we are required to assume that the boundary of $\Omega$, $\partial\Omega$, is sufficiently piecewise smooth, so that the results present in~\cite{MR1014883} guarantee the existence of an interpolate of the true solution of~(\ref{1}) that satisfies the required error bound (see Theorem~\ref{error:est:thm} for further details). Furthermore, in order to define our numerical method, it is necessary to consider the tangential Laplacian (see Lemma~\ref{tange:idents}), on portions of the boundary. This requires the assumption that $\Omega$ is at least piecewise $C^{2,\alpha}$, $\alpha\in(0,1)$.
\end{remark}
In Section~\ref{Error estimates}, we require the following \emph{error estimate necessary mesh assumptions:}
we assume that the family of triangulations $\{\Th\}_h$ is \emph{regular of order $m$,} that is, we assume that conditions~(\ref{meshcond1}) and~(\ref{meshcond2})  hold, and for some integer $m\ge2$, for any $h$, and any element $K\in\Th$ one has $F_K\in C^{m+1}$, and
$$\sup_h\sup_{K\in\Th}\sup_{\hat x\in\hat K}\|\nabla^lF_K(\hat x)\|\|\tilde{B}_K\|^{-l}<\infty,\quad 2\le l\le m+1.$$
Note that since $F_K=\tilde{F}_K+\Phi_K$, where $\tilde{F}_K$ is an affine map, it follows that $\nabla^lF_K=\nabla^l\Phi_K$ for all integers $2\le l\le m+1$.

\emph{Jump and average operators}.
For each face $F=\overline{K}\cap\overline{K'}$ for some $K,K'\in\Th$ (in the case that $F\in\mathscr{E}^b_h$ take $F=\partial K$), with corresponding unit normal vector $n_F$ (which, for convention is chosen so that it is the outward normal to $K$), we define the jump operator, $\jump{\cdot}_F$ over $F$, by
\begin{equation*}
\jump{\phi}_F = \left\{
\begin{aligned}
&(\phi|_K)|_F-(\phi|_{K'})|_F\,\,\mbox{if}\,\,F\in\Ei,\\
&(\phi|_K)|_F\,\,\mbox{if}\,\,F\in\Eb,\\
\end{aligned}
\right.
\end{equation*} 
and the average operator, $\avg{\cdot}$, by
\begin{equation*}
\avg{\phi}_F = \left\{
\begin{aligned}
&\frac{1}{2}((\phi|_K)|_F+(\phi|_{K'})|_F)\,\,\mbox{if}\,\,F\in\Ei,\\
&(\phi|_K)|_F\,\,\mbox{if}\,\,F\in\Eb.\\
\end{aligned}
\right.
\end{equation*} 
For two matrices $A,B\in\mathbb R^{d\times d}$, we set $A:B:=\sum_{i,j=1}^dA_{ij}B_{ij}$. For an element $K$, we define the inner product $\langle\cdot,\cdot\rangle_K$ by
\begin{equation*}
\langle u,v\rangle_K := \left\{
\begin{aligned}
&\int_Ku\,v\,\,\mbox{if}\,\,u,v\in L^2(K),\\
&\int_Ku\cdot v\,\,\mbox{if}\,\,u,v\in L^2(K;\mathbb R^d),\\
&\int_Ku:v\,\,\mbox{if}\,\,u,v\in L^2(K;\mathbb R^{d\times d}).
\end{aligned}
\right.
\end{equation*} 
Any ambiguity in this notation will be resolved by the arguments of the bilinear form. The bilinear forms $\langle\cdot,\cdot\rangle_{\partial K}$ and $\langle\cdot,\cdot\rangle_{F}$ for $F\in\Eib$, are defined similarly. Note that for $F\in\Eb$, in general, to calculate an $\langle\cdot,\cdot\rangle_F$ (with entries that ensure that the integral is well defined), one can utilise the corresponding map $F_K$ (where $F\subset\partial K$), and apply a change of variables, yielding an integral on a face of the reference simplex, $\hat{K}$.

\emph{Function spaces}. For each $K\in\Th$, recall that $\mathbb P^{p}(K)$ is the space of all polynomials with either total or partial degree less than or equal to $p$. The discontinuous Galerkin finite element space $V_{h,p}$ is defined by
\begin{equation}\label{FEspace:def}
V_{h,p}:=\{v\in L^2(\Omega):~v|_K=\rho\circ F_K^{-1},\,\rho\in\mathbb P^{p}(\hat K),~\forall K\in\Th\},
\end{equation}
where $p\in\mathbb N$, $p\ge 2$.
Let $\mathbf s=(s_K:K\in\Th)$ denote a vector of nonnegative real numbers and let $r\in[1,\infty]$.\\
The broken Sobolev space $W^{\mathbf s,r}(\Omega;\Th)$ is defined by
\begin{equation}
W^{\mathbf s,r}(\Omega;\Th):=\{v\in L^2(\Omega):~v|_K\in W^{s_K,r}(K)~\forall K\in\Th\}.
\end{equation}
We denote $H^\mathbf{s}(\Omega;\Th):=W^{\mathbf s,2}(\Omega;\Th)$, and set $W^{s,r}(\Omega;\Th):=W^{\mathbf s,r}(\Omega;\Th)$, in the case that $s_K=s,~s\ge 0$, for all $K\in\Th$. For $v\in W^{1,r}(\Omega;\Th)$, let 
$\nabla_hv\in L^r(\Omega;\mathbb R^d)$ denote the discrete (also known as broken) gradient of $v$, i.e., $(\nabla_hv)|_K=\nabla(v|_K)$ for all $K\in\Th$. Higher order discrete derivatives are defined in a similar way. We define a norm on $W^{s,r}(\Omega;\Th)$ by
\begin{equation}
\|v\|^r_{W^{s,r}(\Omega;\Th)}:=\sum_{K\in\Th}\|v\|^r_{W^{s,r}(K)}
\end{equation}
with the usual modification when $r=\infty$.
\emph{Tangential differential operators.} For $F\in\mathscr{E}^{i,b}$, denote the space of $H^s$-regular tangential vector fields on $F$ by $H^s_\Ta(F):=\{v\in H^s(F)^d:v\cdot n_F=0~\mbox{on}~F\}$. Below we define the tangential gradient 
$\nabla_\Ta:H^s(F)\to H^{s-1}_\Ta(F)$ and the tangential divergence $\operatorname{div}_\Ta:H^s_\Ta(F)\to H^{s-1}(F)$, where $1\le s\le 3$ (note that in the case that $\partial\Omega$ is piecewise $C^{m,\alpha}$, with $m\ge 3$, $\alpha\in(0,1)$, we are able to consider $1\le s\le m$). Let $\{t_i\}_{i=1}^{d-1}\subset\mathbb R^d$ be an orthonormal coordinate system on $F$. Then, for $u\in H^s(F)$ and $v=\sum_{i=1}^{d-1}v_it_i$, with $v_i\in H^s(F)$ for $i=1,\ldots,d-1$, we define
\begin{equation}
\nabla_\Ta u:=\sum_{i=1}^{d-1}t_i\frac{\partial u}{\partial t_i},\quad\operatorname{div}_\Ta v:=\sum_{i=1}^{d-1}\frac{\partial v_i}{\partial t_i}.
\end{equation}

\begin{definition}[Mean curvature]
Let $F$ be a $C^{2,\alpha}$, portion of $\partial\Omega$, with $\alpha\in(0,1)$. Then, $F$ defines a $(d-1)$--dimensional hypersurface in $\mathbb R^d$. We define the mean curvature, $\mathcal{H}_F$ of $F$ as follows
\begin{equation*}
\mathcal{H}_F = (\nabla\cdot n_{F})|_F= (\operatorname{Tr}(\nabla n_{F}^T))|_F,
\end{equation*}
where $n_{F}$ is the unit outward normal to $F$.
\end{definition}
\begin{lemma}\label{tange:idents}
Let $\Omega$ be a piecewise $C^{2,\alpha}$-piecewise convex domain (in the sense of Definition~\ref{p-conv-def}), and let $\mathscr{T}_h$ be a mesh on $\Omega$ consisting of \emph{possibly curved} simplices or parallelepipeds. Then, for each $K\in\mathscr{T}_h$ and each face $F\subset\partial K$, the following identities hold:
\begin{equation}\label{identities}
\begin{aligned}
\tau_F(\nabla v)  &= \nabla_\Ta(\tau_Fv)+\left(\tau_F\frac{\partial v}{\partial n_F}\right)n_F\quad\forall v\in H^s(K),\,s>3/2,\\
\tau_F(\Delta v)  &= \operatorname{div}_\Ta\nabla_\Ta(\tau_Fv)+\mc_F\left(\tau_F\frac{\partial v}{\partial n_F}\right)+\tau_F\frac{\partial}{\partial n_F}(\nabla v\cdot n_F),\quad\forall v\in H^s(K),\,s>5/2,
\end{aligned}
\end{equation}
where $\mc_F$ is the mean curvature of the face $F$, and $\tau_F$ is the trace operator associated to $F$.
\end{lemma}
\emph{Proof:} First, if $F\in\Ei$, then $F$ is flat, and both identities in~(\ref{identities}) follow from Lemma 4 in~\cite{MR3077903}. 

If $F\in\Eb$, then the identities follow similarly to the proof of Lemma 4 in~\cite{MR3077903}, taking into account the fact that for $K\in\Th$ such that $\Eb\ni F\subset\partial K$, the Laplacian of a smooth function $u\in C^\infty(\overline{K})$ can be decomposed as follows
$$\Delta u|_F =\operatorname{div}_\Ta\nabla_\Ta u+\left.\mc_F\frac{\partial u}{\partial n}\right|_F+\left.\frac{\partial^2u}{\partial n^2}\right|_F.$$
Noting that the trace operator, $\tau$, commutes with partial derivatives, by the density of smooth functions in $H^s(K)$, we obtain
$$\tau_F(\Delta v)  = \operatorname{div}_\Ta\nabla_T(\tau_Fv)+\mc_F\left(\tau_F\frac{\partial v}{\partial n_F}\right)+\tau_F\frac{\partial}{\partial n_F}(\nabla v\cdot n_F).\quad\quad\square$$
\subsection{Trace inverse estimate}
We prove the following result, based on the proof of Lemma 2.4 from~\cite{MR1014883}.
\begin{lemma}
Assume that the simplex $K=F_K(\hat K)$, where $F_K\in C^k$, $k\in\mathbb N$, and that $F\in\Eib$ is a face of $K$. Then, for any $v\in V_{h,p}$, the trace of $v$, $\tau_F(v|_K)\in C^k(\overline{F})$, and we have, for any integer $0\le m\le k$,
\begin{equation}\label{bnd:inv:est}
\|v\|_{H^m(F)}\le C\tilde{h}_F^{-m}\|v\|_{L^2(F)},
\end{equation}
where $C$ is a positive constant independent of the mesh size $h$.
\end{lemma}
\emph{Proof:} From the second inequality in the proof of Lemma 2.4 in~\cite{MR1014883}, for $v\in V_{h,p}$, we obtain the following
\begin{equation*}
|v|^2_{H^m(F)}\le C|\deter\tilde{B}_F|\left(\sum_{r=1}^m\|\tilde{B}_K\|^{4(m-r)}\|\tilde{B}_K^{-1}\|^{2m}|\hat v|_{H^r(\hat F)}^2\right),
\end{equation*}
where $\tilde{B}_F$ is the restriction of $\tilde{B}_K$ to $\hat F:=F_K^{-1}(F)$, and $\hat v = v\circ F_K^{-1}$.

As all norms are equivalent in finite dimensions, we see that, for $1\le r\le m$,
$$|\hat v|_{H^r(\hat F)}\le \|\hat v\|_{H^r(\hat F)}\le C\|\hat v\|_{L^2(\hat F)},$$
and so
\begin{equation}\label{nearly:done}
\begin{split}
|v|^2_{H^m(F)}&\le C|\deter\tilde{B}_F|\left(\sum_{r=1}^m\|\tilde{B}_K\|^{4(m-r)}\right)\|\tilde{B}_K^{-1}\|^{2m}|\hat v|_{L^2(\hat F)}^2\\
&\le C|\deter\tilde{B}_F|\|\tilde{B}_K^{-1}\|^{2m}\|\hat v\|_{L^2(\hat F)}^2.
\end{split}
\end{equation}
Now, note that
\begin{equation}\label{nearly:done-1}
\|\hat v\|_{L^2(\hat F)}^2 = \int_{\hat F} (v\circ F_K)^2\le \sup_{x\in F}\left|\deter DF^{-1}_K(x)\right|\left(\int_{\hat F} v^2\right) =\sup_{x\in F}\left|\deter DF^{-1}_K(x)\right|\|v\|_{L^2(F)}^2.
\end{equation}
Furthermore,
$$DF_K^{-1} = \tilde{B}_K^{-1}(I_d+D\Phi_K\tilde{B}_K^{-1})^{-1},$$
where the inverse of $I_d+D\Phi_K\tilde{B}_K^{-1}$ exists, due to the fact that
\begin{equation}\label{||A||:bound}
\sup_{\hat x\in\hat K}\|D\Phi_K(\hat x)\tilde{B}_K^{-1}\|\le C_K<1.
\end{equation}
From this, we obtain
$$\left|\deter DF_K^{-1}\right| = \left|\deter\tilde{B}_K^{-1}\right|\left|\deter(I_d+D\Phi_K\tilde{B}_K^{-1})\right|^{-1}.$$
Then, one can see that~(\ref{||A||:bound}) implies that
$$\deter(I_d+D\Phi_K\tilde{B}_K^{-1})|_{\hat F}\ge(1-C_K)^{d-1},$$ 
which yields
$$\sup_{x\in F}|\deter DF_K^{-1}|\le\frac{|\deter\tilde{B}_F^{-1}|}{(1-C_K)^{d-1}}\le C|\deter\tilde{B}_F^{-1}|=C|\deter\tilde{B}_F|^{-1}$$
Applying this estimate, in conjunction with~(\ref{nearly:done-1})  to~(\ref{nearly:done}), gives us the following:
$$|v|^2_{H^m(F)}\le C\|\tilde{B}_K^{-1}\|^{2m}|v|_{L^2(F)}^2\le C\tilde{h}_F^{-2m}|v|_{L^2(F)}^2.$$
Taking square roots in the estimate above, we obtain the desired result.$\quad\quad\square$

\subsection{Numerical scheme} The definition of the numerical scheme requires the following bilinear forms. Firstly, the bilinear form $B_{h,*}:V_{h,p}\times V_{h,p}\to\mathbb R$ is defined by
\begin{equation}\label{Bh*def}
\begin{aligned}
B_{h,*}(u_h,v_h):=&\sum_{K\in\Th}\langle D^2u_h,D^2v_h\rangle_K\\
&+\sum_{F\in\mathscr{E}^{i}_h}[\langle \operatorname{div}_\Ta\nabla_\Ta\average{u_h},\jump{\nabla v_h\cdot n_F}\rangle_F+\langle \operatorname{div}_\Ta\nabla_\Ta\average{v_h},\jump{\nabla u_h\cdot n_F}\rangle_F]\\
&-\sum_{F\in\mathscr{E}^{i,b}_h}[\langle\nabla_\Ta\average{\nabla u_h\cdot n_F},\jump{\nabla_\Ta v_h}\rangle_F+\langle\nabla_\Ta\average{\nabla v_h\cdot n_F},\jump{\nabla_\Ta u_h}\rangle_F]\\
&+\sum_{F\in\mathscr{E}^b_h}\left\langle\mc_F\frac{\partial u_h}{\partial n_F},\frac{\partial v_h}{\partial n_F}\right\rangle_F+\langle\nabla_\Ta u_h,\nabla n^T_F\nabla_\Ta v_h\rangle_F+\left\langle\frac{\partial u_h}{\partial n_F},n_F^T\nabla n^T_F\nabla_\Ta v_h\right\rangle_F,
\end{aligned}
\end{equation}
where $\mc_F$ is the mean curvature of the face $F$, and $u_h,~v_h\in V_{h,p}$ throughout this section. Then, for positive edge-dependent quantities $\mu_F$ and $\eta_F$ to be specified later, the jump stabilization bilinear form $J_h:V_{h,p}\times V_{h,p}\to\mathbb R$ is defined by 
\begin{equation}
\begin{aligned}
J_h(u_h,v_h):=&\sum_{F\in\mathscr{E}^{i,b}_h}[\mu_F\langle\jump{\nabla_\Ta u_h},\jump{\nabla_\Ta v_h}\rangle_F+\eta_F\langle\jump{u_h},\jump{v_h}\rangle_F]\\
&+\sum_{F\in\mathscr{E}^{i}_h}\mu_F\langle\jump{\nabla u_h\cdot n_F},\jump{\nabla v_h\cdot n_F}\rangle_F.
\end{aligned}
\end{equation}
For each $\theta\in(0,1]$, we define the bilinear form $B_{h,\theta}:V_{h,p}\times V_{h,p}\to\mathbb R$ by
\begin{equation}
B_{h,\theta}(u_h,v_h):=\theta B_{h,*}(u_h,v_h)+(1-\theta)\sum_{K\in\Th}\langle\Delta u_h,\Delta v_h\rangle_K+J_h(u_h,v_h).
\end{equation}
Finally, the bilinear form $A_h:V_{h,p}\times V_{h,p}\to\mathbb R$ is defined by
\begin{equation}\label{A:def}
A_h(u_h,v_h):=\sum_{K\in\Th}\langle\gamma A:D^2u_h,\Delta v_h\rangle_K+B_{h,1/2}(u_h,v_h)-\sum_{K\in\Th}\langle\Delta u_h,\Delta v_h\rangle_K.
\end{equation}
The scheme for approximating the solution of~(\ref{1}) is to find $u_h\in V_{h,p}$ such that
\begin{equation}\label{method}
A_h(u_h,v_h)=\sum_{K\in\Th}\langle\gamma f,\Delta v_h\rangle_K\quad\forall w_h\in V_{h,p}.
\end{equation}
\begin{remark}
In~(\ref{Bh*def})-(\ref{A:def}) we have defined the bilinear forms $B_{h,*},J_h,B_{h,\theta},A_h:V_{h,p}\times V_{h,p}\to\Reals$. The main difference between these bilinear forms, and the bilinear forms given presented in Section 3 of~\cite{MR3077903} is in the bilinear form $B_{h,*}$ (and thus, by definition, in $B_{h,\theta}$ and $A_h$). 
In particular, the bilinear form $B_{h,*}$ (given by~(\ref{Bh*def})) contains the following extra terms:
$$\sum_{F\in\mathscr{E}^b_h}\left\langle\mc_F\frac{\partial u_h}{\partial n_F},\frac{\partial v_h}{\partial n_F}\right\rangle_F+\langle\nabla_\Ta u_h,\nabla n^T_F\nabla_\Ta v_h\rangle_F+\left\langle\frac{\partial u_h}{\partial n_F},n_F^T\nabla n^T_F\nabla_\Ta v_h\right\rangle_F,$$
which arise due to the curvature of $\partial\Omega$. Indeed, if $\Omega$ is polytopal (which is a necessary assumption of Lemmas 5, 7, and 8, as well as Theorems 8 and 9 in~\cite{MR3077903}), then all of the faces $F\in\Eb$ are flat, and so 
$$\mc_F = 0\quad\mbox{and}\quad\nabla n_F = 0,\quad\mbox{for all }F\in\Eb,$$
which means that the additional terms vanish. In experiment~\ref{exp:cons}, the results imply the necessity of these extra terms when $\partial\Omega$ has curved boundary portions. Furthermore, the presence of these additional terms requires the application of new techniques, in order to prove that the numerical method admits a unique solution (see Theorem~\ref{mainthm}), is consistent (see Lemma~\ref{cons:lemma}), and that the resulting solution satisfies optimal error bounds (see Theorem~\ref{error:est:thm}).
\end{remark}
\begin{remark} For particular geometries, the final boundary terms of~(\ref{Bh*def}) may simplify. For example, if $\nabla n_F^T$ is symmetric for some $F\in\Eb$, we see that
$$\left\langle\frac{\partial u_h}{\partial n_F},n_F^T\nabla n^T_F\nabla_\Ta v_h\right\rangle_F=\left\langle\gradTa v_h,(\nabla n_F^Tn_F)\frac{\partial u_h}{\partial n_F}\right\rangle_F = 0.$$
\end{remark}
\begin{example}[Sphere of radius $r$]
In the case that $\Omega = \{ x\in\mathbb R^d:|x|<r\}$, we see that $$\partial\Omega = \{ x\in\mathbb R^d:|x|=r\}=r\mathbb S^{d-1}.$$
Here, the unit outward normal to $\partial\Omega$ is given by $n=x/|x|$, and thus
$$D_in_j = \frac{r^2\delta_{ij}-x_ix_j}{r^3}=D_jn_i,\quad i,j=1,\ldots,d,$$
where $\delta_{ij}$ is the Kronecker-delta symbol.
From this, we can calculate the following:
\begin{equation*}
\begin{split}
\mc & = \frac{d-1}{r},\\
(\gradTa u_h)^T\nabla n^T\gradTa v_h &= \frac{d-1}{r}(\gradTa u_h)^T\gradTa v_h,\\
\end{split}
\end{equation*}
and so
\begin{equation*}
\sum_{F\in\mathscr{E}^b_h}\left\langle\mathcal{H}\frac{\partial u_h}{\partial n_F},\frac{\partial v_h}{\partial n_F}\right\rangle_F+\langle\nabla_\Ta u_h,\nabla n^T_F\nabla_\Ta v_h\rangle_F+\left\langle\frac{\partial u_h}{\partial n_F},n_F^T\nabla n^T_F\nabla_\Ta v_h\right\rangle_F = 
\frac{d-1}{r}\sum_{F\in\mathscr{E}^b_h}\langle\nabla u_h,\nabla v_h\rangle_F.
\end{equation*}
\end{example}
\end{section}
\begin{section}{Analysis of the numerical method}\label{Analysis of the numerical method}
We will first prove that the method is stable, yielding a result for the existence and uniqueness of a numerical solution to~(\ref{method})
\subsection{Stability of the method}
Let $c_{*}$ be a positive constant independent of $h$ and to be determined later. For each $\theta\in(0,1]$ define the functional $\|\cdot\|_{h,\theta}:V_{h,p}\to\mathbb R^+$ by
\begin{equation}
\|v_h\|_{h,\theta}^2:=\sum_{K\in\Th}[\theta|v_h|_{H^2(K)}^2+(1-\theta)\|\Delta v_h\|^2_{L^2(K)}]+\frac{\theta}{2}\sum_{F\in\Eb}\|\mc_F^{1/2}\frac{\partial v_h}{\partial n}\|_{L^2(F)}^2+c_{*}J_h(v_h,v_h).
\end{equation}
\begin{lemma}
For any $\theta\in(0,1]$, $\|\cdot\|_{h,\theta}$ is a norm on $V_{h,p}$.
\end{lemma}
\emph{Proof:} Homogeneity and the triangle inequality are clear. It remains to show that if $\|v_h\|_{h,\theta}=0$, then $v_h=0$ for $v_h\in V_{h,p}$. Let $v_h\in V_{h,p}$ satisfy $\|v_h\|_{h,\theta}=0$ for some $\theta\in(0,1]$. Since $\theta\in(0,1]$, it follows that
$|v_h|_{H^2(\Omega;\Th)} = 0$, and thus $v_h$ must be piecewise affine. Furthermore, $J_h(v_h,v_h)=0$ implies that
$\jump{\nabla v_h}= 0$ for all $F\in\mathscr{E}^{i}_h$, and $\jump{v_h}= 0$ for all $F\in\mathscr{E}^{i,b}_h$.
It follows that $v_h$ is an affine function that satisfies $v_h|_{\partial\Omega} = 0$, and so $v_h\equiv 0$.$\quad\quad\square$

The following stability lemma requires some assumptions upon the piecewise nature of $\partial\Omega$, namely that $$\partial \Omega = \bigcup_{n=1}^N\overline{\Gamma}_n,$$
for some $\Gamma_1,\ldots,\Gamma_N$, where each $\Gamma_n$ admits a local representation as the graph of $C^{2,\alpha}$ function. Moreover, for each $n$, the unit outward normal to $\Gamma_n$, $n_{\Gamma_n}$, must satisfy
\begin{equation}\label{star-14/12}
\nabla n_{\Gamma_n}(x) = 0\,\,\forall x\in\Gamma_n\quad\mbox{or}\quad\nabla n_{\Gamma_n}(x) \ne0\,\,\forall x\in\Gamma_n.
\end{equation}
Coupling this with the assumption that $\Omega$ is piecewise $C^{2,\alpha}$-piecewise convex, $\alpha\in(0,1)$, leads to each portion $\Gamma_n$ of $\partial\Omega$ to either be flat, or curved, with nonvanishing \emph{positive} curvature. This means that our framework encompasses polytopal domains, curved domains, and domains with both curved boundary portions and flat boundary portions (for example, the key-hole shaped domain~(\ref{Key-hole})).
\begin{lemma}\label{lemma6}
Let $\Omega$ be a piecewise $C^{2,\alpha}$-piecewise convex domain, with $\alpha\in(0,1)$, satisfying~(\ref{star-14/12}), and let $\{\Th\}_h$ be a regular sequence of (possibly curved) simplicial or parallelepipedal meshes satisfying~(\ref{meshcond1})-(\ref{Fcontained}). Then, for each constant $\kappa>1$, there exists a positive constant $c_{\operatorname{stab}}$, independent of $h$, $p$, and $\theta$, such that for any $v_h\in V_{h,p}$ and any $\theta\in(0,1]$, we have
\begin{equation}\label{state:lemma}
\kappa B_{h,\theta}(v_h,v_h)\ge\theta|v_h|^2_{H^2(\Omega;\Th)}+(1-\theta)\sum_{K\in\Th}\|\Delta v_h\|^2_{L^2(K)}+\frac{1}{2}J_h(v_h,v_h)+\frac{\theta}{2}\sum_{F\in\Eb}\|\mc_F^{1/2}\frac{\partial v_h}{\partial n}\|_{L^2(F)}^2,
\end{equation}
where, for some fixed constant $\sigma\ge 1$, the jump penalty parameters $\mu_F$ and $\eta_F$ satisfy
\begin{equation}\label{require}
\mu_F=\sigma \left(c_{\operatorname{stab}}\frac{1}{\tilde{h}_F}+c_\mathcal{H}\right)\quad\mbox{and}\quad\eta_F\ge\frac{\sigma}{\tilde{h}_F^3}.
\end{equation}
Here $c_\mathcal{H}$ depends on the mean curvature lower bound, which we define by
$$\mathcal{H}_{\operatorname{min}}:=\min_{F\in\Eb:\nabla n_F\ne 0}\inf_F\mc_F,
$$
 and the upper bound on the tangential gradient of the normal vector, $\max_{F\in\Eb}\|\nabla_\Ta n_F^T\|_{L^\infty(F)}$.
\end{lemma}
\emph{Proof:} The proof is similar to that of~\cite{MR3196952}, Section 6, Lemma 6; in this case we must now deal with the extra terms arising in the bilinear form $B_{h,*}$ due to the curvature of the boundary, $\partial\Omega$.

Firstly, for $v_h\in V_{h,p}$, we have
$$B_{h,\theta}(v_h,v_h) = \theta|v_h|^2_{H^2(\Omega;\mathscr{T}_h)}+(1-\theta)\sum_{K\in\mathscr{T}_h}\|\Delta v_h\|_{L^2(K)}^2+J_h(v_h,v_h)+\theta\sum_{i=1}^6I_i,$$
where
\begin{equation*}
\begin{split}
I_1&:=2\sum_{F\in\mathscr{E}^i_h}\langle\operatorname{div}_\Ta \nabla_\Ta \average{v_h},\jump{\nabla v_h\cdot n_F}\rangle_F,\quad I_2:=\sum_{F\in\mathscr{E}^b_h}\left\langle\mc_F\frac{\partial v_h}{\partial n_F},\frac{\partial v_h}{\partial n_F}\right\rangle_F,\\
I_3&:=-2\sum_{F\in\mathscr{E}^{i}_h}\langle\nabla_\Ta \average{\nabla v_h\cdot n_F},\jump{\nabla_\Ta v_h}\rangle_F,\quad I_4:=-2\sum_{F\in\mathscr{E}^{b}_h}\langle\nabla_\Ta \average{\nabla v_h\cdot n_F},\jump{\nabla_\Ta v_h}\rangle_F,\\
\quad I_5&:=\sum_{F\in\mathscr{E}^b_h}\langle\nabla_\Ta v_h,\nabla n^T_F\nabla_\Ta v_h\rangle_F,\quad
I_6:=\sum_{F\in\mathscr{E}^b_h}\left\langle\frac{\partial v_h}{\partial n_F},n^T_F\nabla n^T_F\nabla_\Ta v_h\right\rangle_F. 
\end{split}
\end{equation*}
In~\cite{MR3077903}, it is shown that there is a constant $C(d)$ depending only on $d$, such that, for any $\delta>0$,
\begin{equation}
\begin{split}
&|I_1|\le\delta C(d)C_{\operatorname{Tr}}c_\mathcal{F}\sum_{K\in\Th}\|D^2v_h\|^2_{L^2(K)}+\sum_{F\in\mathscr{E}^i_h}\frac{1}{\delta\tilde{h}_F}\|\jump{\nabla v_h\cdot n_F}\|^2_{L^2(F)},\\
&|I_3|\le\delta C(d)C_{\operatorname{Tr}}c_\mathcal{F}\sum_{K\in\Th}\|D^2v_h\|^2_{L^2(K)}+\sum_{F\in\mathscr{E}^{i}_h}\frac{1}{\delta\tilde{h}_F}\|\jump{\nabla_\Ta v_h}\|^2_{L^2(F)},
\end{split}
\end{equation}
where $C_{\operatorname{Tr}}$ is the combined constant of the trace and inverse inequalities, and $c_\mathcal{F}$ is given by~(\ref{meshcond1}). We shall prove a similar bound for $I_4$ by noting that, for any $F\in\Eb$,
\begin{equation*}
\begin{split}
\gradTa\left(\frac{\partial v_h}{\partial n_F}\right) & = \sum_{k=1}^{d-1}\frac{\partial}{\partial t_k}\left(\frac{\partial v_h}{\partial n_F}\right)t_k\\
& = \sum_{k=1}^{d-1}(t_k)^T\nabla\left(\frac{\partial v_h}{\partial n_F}\right)t_k\\
& = \sum_{k=1}^{d-1}((t_k)^TD^2v_h\,n_F+(t_k)^T\nabla n_F^T\nabla v_h)t_k\\
& = \sum_{k=1}^{d-1}((t_k)^TD^2v_h\,n_F+(t_k)^T\nabla n_F^T(n_F\frac{\partial v_h}{\partial n_F}+\gradTa v_h))t_k\\
& = \sum_{k=1}^{d-1}((t_k)^TD^2v_h\,n_F+(t_k)^T\nabla n_F^T\gradTa v_h)t_k,\\
\end{split}
\end{equation*}
and so
\begin{equation}\label{gradTdn:bound}
|\gradTa\left(\frac{\partial v_h}{\partial n_F}\right)|\le(d-1)(|D^2v_h|+|\nabla n_F||\gradTa v_h|).
\end{equation}
Now, we see that
\begin{equation*}
\begin{split}
|I_4|& = 2\left|\sum_{F\in\Eb}\left\langle\gradTa\left(\frac{\partial v_h}{\partial n_F}\right),\gradTa v_h\right\rangle_F\right|\\
&\le 2\sum_{F\in\Eb}\left\|\gradTa\left(\frac{\partial v_h}{\partial n_F}\right)\right\|_{L^2(F)}\|\gradTa v_h\|_{L^2(F)}\\
&\le2(d-1)\sum_{F\in\Eb}\left[\|D^2v_h\|_{L^2(F)}\|\gradTa v_h\|_{L^2(F)}+\max_{F\in\Eb}\|\gradTa n_F^T\|_{L^\infty(F)}\|\gradTa v_h\|_{L^2(F)}^2\right]\\
&\le(d-1)\sum_{F\in\Eb}\left[\tilde{h}_F\|D^2v_h\|_{L^2(F)}^2+(2\max_{F\in\Eb}\|\gradTa n_F^T\|_{L^\infty(F)}+\tilde{h}_F^{-1})\|\gradTa v_h\|_{L^2(F)}^2\right].
\end{split}
\end{equation*}
Then, applying~(\ref{bnd:inv:est}) with $m = 2$, we obtain
$$|I_4|\le C\sum_{F\in\Eb}\left[\tilde{h}_F^{-3}\|v_h\|_{L^2(F)}^2+(\max_{F\in\Eb}\|\gradTa n_F^T\|_{L^\infty(F)}+\tilde{h}_F^{-1})\|\gradTa v_h\|_{L^2(F)}^2\right].$$
One can also see that
$$I_2=\sum_{F\in\Eb}\left\|\mc_F^{1/2}\frac{\partial v_h}{\partial n}\right\|_{L^2(F)}^2.$$
To see that the value $I_5$ is nonnegative, we first note that by mesh assumption~(\ref{Fcontained}), each $F\in\Eb$ is contained in $\Gamma_n$, for some $n\in\{1,\ldots,N\}$. Thus, the unit normal to $F$, $n_F$, corresponds to the unit normal to $\Gamma_n$, $n_{\Gamma_n}$, and so, it follows that
$$n_F = n_{\Gamma_n}= \frac{\nabla\varphi_n}{|\nabla\varphi_n|},$$
due to~(\ref{convex-normal}), for a $C^{2,\alpha}$ convex function $\varphi_n$, $\alpha\in(0,1)$. From this we can calculate:
\begin{equation*}
\begin{split}
[\nabla n_F^T]^i_j = D_i\left(\frac{D_j\varphi_n}{|\nabla\varphi_n|}\right) & = \frac{|\nabla\varphi_n|\,D^2_{ij}\varphi_n-D_j\varphi_n\,D_i((\sum_{k=1}^d(D_k\varphi_n)^2)^{1/2})}{|\nabla\varphi_n|^2}\\
& =\frac{|\nabla\varphi_n|\,D^2_{ij}\varphi_n-D_j\varphi_n\sum_{k=1}^dD_k\varphi_n\, D^2_{ik}\varphi_n/|\nabla\varphi_n|}{|\nabla\varphi_n|^2}\\
& = \frac{|\nabla\varphi_n|^2\,D^2_{ij}\varphi_n-D_j\varphi_n\sum_{k=1}^dD_k\varphi_n\,D^2_{ik}\varphi_n}{|\nabla\varphi_n|^3}.
\end{split}
\end{equation*}
Now let $\tau,\xi$ be two tangent vectors to $\partial\Omega$; we then see that
\begin{equation*}
\begin{split}
\tau^T\nabla n_F^T\,\xi & = \frac{|\nabla\varphi_n|^2\tau^TD^2\varphi_n\,\xi-\sum_{i,j,k=1}^dD^2_{ik}\varphi_n\,D_k\varphi_n\,D_j\varphi_n\,\tau_i\xi_j}{|\nabla\varphi_n|^3}\\
& =  \frac{|\nabla\varphi_n|^2\tau^TD^2\varphi_n\,\xi}{|\nabla\varphi_n|^3}-\frac{\sum_{i,k=1}^dD^2_{ik}\varphi_n\,D_k\varphi_n\,\tau_i\sum_{j=1}^d\frac{D_j\varphi_n\,\xi_j}{|\nabla\varphi_n|}}{|\nabla\varphi_n|^2}\\
& = \frac{|\nabla\varphi_n|^2\tau^TD^2\varphi_n\,\xi}{|\nabla\varphi_n|^3}-\frac{\sum_{i,k=1}^dD^2_{ik}\varphi_n\,D_k\varphi_n\,\tau_i(n\cdot\xi)}{|\nabla\varphi_n|^2}\\
& = \frac{|\nabla\varphi_n|^2\tau^TD^2\varphi_n\,\xi}{|\nabla\varphi_n|^3}.
\end{split}
\end{equation*}
Recall that the function $\varphi_n$ is convex, and so its Hessian is positive semidefinite, noting the above calculation, and taking into account the fact that the tangential gradient, $\nabla_\Ta $, of a smooth function $w$ is a tangent vector-valued function, we obtain, for any face $F\in\mathscr{E}^b_h,$
$$(\gradTa w)^T\nabla n_F^T\gradTa w=\frac{(\gradTa w)^TD^2\varphi_n\gradTa w}{|\nabla\varphi_n|}\ge0\,\,\mbox{on}\,\,F,\mbox{ for some }n\in\{1,\ldots,N\}.$$
This inequality extends to $v_h\in V_{h,p}$ by construction of the trace operator, and thus we find that
\begin{equation}\label{I5}
I_5 = \sum_{F\in\mathscr{E}^b_h}\langle \gradTa v_h,\nabla n_F^T\gradTa v_h\rangle_F\ge 0.
\end{equation}
Finally, for $I_6$, we use the Cauchy--Schwarz inequality with a parameter, to obtain
\begin{equation*}
\begin{aligned}
I_6& = \sum_{F\in\Eb}\left\langle\frac{\partial v_h}{\partial n_F},n_F^T\nabla n_F^T\nabla_\Ta v_h\right\rangle_F\\
& = \sum_{F\in\Eb:\nabla n_F\ne 0}\left\langle\frac{\partial v_h}{\partial n_F},n_F^T\nabla n_F^T\nabla_\Ta v_h\right\rangle_F\\
& = \sum_{F\in\Eb:\nabla n_F\ne 0}\left\langle\mc_F^{1/2}\frac{\partial v_h}{\partial n_F},\mc_F^{-1/2}n_F^T\nabla n_F^T\nabla_\Ta v_h\right\rangle_F\\
& \ge-\sum_{F\in\mathscr{E}^b_h:\nabla n_F\ne 0}\left(\frac{1}{2}\|\mc_F^{1/2}\frac{\partial v_h}{\partial n_F}\|^2_{L^2(F)}+\frac{\sup_{F\in\mathscr{E}_h^b}\|\nabla n^T_F\|_{L^\infty(F)}}{2\mathcal{H}_{\operatorname{min}}}\|\gradTa v_h\|^2_{L^2(F)}\right)\\
& \ge-\sum_{F\in\mathscr{E}^b_h}\left(\frac{1}{2}\|\mc_F^{1/2}\frac{\partial v_h}{\partial n_F}\|^2_{L^2(F)}+\frac{\sup_{F\in\mathscr{E}_h^b}\|\nabla n^T_F\|_{L^\infty(F)}}{2\mathcal{H}_{\operatorname{min}}}\|\gradTa v_h\|^2_{L^2(F)}\right).
\end{aligned}
\end{equation*}
Now that we have lower bounds on $I_1,\ldots,I_6$, we obtain the following:
\begin{equation*}
B_{h,\theta}(v_h,v_h)\ge\sum_{i=1}^8A_i,
\end{equation*}
where
\begin{equation*}
\begin{split}
A_1 & = \theta(1-2\delta C(d)C_{\operatorname{Tr}}c_{\mathcal{F}})|v_h|^2_{H^2(\Omega;\Th)},\quad
A_2  = (1-\theta)\sum_{K\in\Th}\|\Delta v_h\|^2_{L^2(K)},\\
A_3 & = \sum_{F\in\Ei}\left(\mu_F-\frac{2\theta}{\delta\tilde{h}_F}\right)\|\jump{\nabla v_h\cdot n_F}\|^2_{L^2(F)},\quad
A_4  = \sum_{F\in\Ei}\left(\mu_F-\frac{\theta}{\delta\tilde{h}_F}\right)\|\jump{\gradTa v_h}\|^2_{L^2(F)},\\
A_5 & = \sum_{F\in\mathscr{E}^b_h}\left(\mu_F-\frac{\theta}{\tilde{h}_F}-\max_{F\in\mathscr{E}_h^b}\|\nabla n^T_F\|_{L^\infty(F)}(1+\frac{1}{2\mathcal{H}_{\operatorname{min}}})\right)\|\gradTa v_h\|_{L^2(F)}^2,\\
A_6 & = \frac{\theta}{2}\sum_{F\in\mathscr{E}^b_h}\|\mc_F^{1/2}\frac{\partial v_h}{\partial n_F}\|^2_{L^2(F)},\quad A_7 = \sum_{F\in\mathscr{E}^{i}_h}\eta_F\|\jump{v_h}\|^2_{L^2(F)},\\
A_8 &= \sum_{F\in\Eb}(\eta_F-\frac{\theta C}{\tilde{h}_F^3})\|v_h\|_{L^2(F)}^2.
\end{split}
\end{equation*}
For any given $\kappa>1$, there is a $\delta>0$ such that $1-2\delta C(d)C_{\operatorname{Tr}}c_{\mathcal{F}}>\kappa^{-1}$. Set $c_{\operatorname{stab}}=2/\delta$, $c_*=\kappa/2$ so that the following inequalities hold for any $\theta\in(0,1]$:
\begin{equation*}
\begin{split}
A_3&\ge\frac{1}{2}\sum_{F\in\Ei}\mu_F\|\jump{\nabla v_h\cdot n_F}\|_{L^2(F)}^2=\kappa^{-1}c_*\sum_{F\in\Ei}\mu_F\|\jump{\nabla v_h\cdot n_F}\|_{L^2(F)}^2,\\
A_4&\ge\frac{1}{2}\sum_{F\in\Ei}\mu_F\|\jump{\gradTa v_h}\|_{L^2(F)}^2=\kappa^{-1}c_*\sum_{F\in\Ei}\mu_F\|\jump{\gradTa v_h}\|_{L^2(F)}^2,\\
A_5&\ge\frac{1}{2}\sum_{F\in\Eb}\mu_F\|\gradTa v_h\|_{L^2(F)}^2=\kappa^{-1}c_*\sum_{F\in\Eb}\mu_F\|\gradTa v_h\|_{L^2(F)}^2,\\
A_7&\ge\frac{1}{2}A_7 = \kappa^{-1}c_*\sum_{F\in\mathscr{E}^{i,b}_h}\eta_F\|\jump{v_h}\|^2_{L^2(F)},\\
A_8&\ge\frac{1}{2}\sum_{F\in\Eb}\eta_F\|v_h\|_{L^2(F)}^2=\kappa^{-1}c_*\sum_{F\in\Eb}\eta_F\|v_h\|_{L^2(F)}^2,
\end{split}
\end{equation*}
whenever $\mu_F$ and $\eta_F$ satisfy~(\ref{require}). Thus we obtain the following
$$\kappa B_{h,\theta}(v_h,v_h)\ge\theta|v_h|^2_{H^2(\Omega;\Th)}+(1-\theta)\sum_{K\in\Th}\|\Delta v_h\|^2_{L^2(K)}+\frac{1}{2}J_h(v_h,v_h)+\frac{\theta}{2}\sum_{F\in\Eb}\|\mc_F^{1/2}\frac{\partial v_h}{\partial n}\|_{L^2(F)}^2.
\quad\quad\square$$
\begin{theorem}\label{mainthm}
Under the hypotheses of Lemma~\ref{lemma6}, let $c_{\operatorname{stab}}$, $c_{\mathcal{H}}$, $\eta_F$ and $\mu_F$ be chosen so that Lemma~\ref{state:lemma} holds with $\kappa<(1-\varepsilon)^{-1}$. Then, for every $v_h\in V_{h,p}$, we have
\begin{equation}
\|v_h\|_{h,1}^2\le\frac{2\kappa}{1-\kappa(1-\epsilon)}A_h(v_h,v_h).
\end{equation}
Therefore, there exists a unique solution $u_h\in V_{h,p}$ of the numerical scheme~(\ref{method}). Furthermore, we have the bound
\begin{equation}
\|u_h\|_{h,1}\le\frac{2\kappa\sqrt{d}\|\gamma\|_{L^\infty(\Omega)}}{1-\kappa^2(1-\varepsilon)}\|f\|_{L^2(\Omega)}.
\end{equation}
\end{theorem}
\emph{Proof:} The proof is the same as the proof of Theorem 8, in~\cite{MR3077903}, Section 4, which relies upon the stability estimate~(\ref{state:lemma}).$\quad\quad\square$

We will now prove a consistency result for our method. This method is central to the error analysis discussed in Section~\ref{Error estimates}, as it allows for a ``Galerkin orthogonality" type argument.
\subsection{Consistency of the method}\label{cons:sec}
\begin{lemma}\label{cons:lemma}
Let $\Omega$ be a piecewise $C^{2,\alpha}$ domain, with $\alpha\in(0,1)$, or a convex polytopal domain, and let $\Th$ be an exact mesh on $\Omega$ consisting of simplices or parallelepipeds possibly with curved boundary faces. Let $w\in H^s(\Omega;\Th)\cap H^2(\Omega)\cap H^1_0(\Omega),$ $s>5/2$.
Then, for every $v_h\in V_{h,p}$, we have the identities
\begin{equation}\label{cons:statement}
B_{h,*}(w,v_h)=\sum_{K\in\Th}\langle\Delta w,\Delta v_h\rangle_K\quad\mbox{and}\quad J_h(w,v_h)=0.
\end{equation}
\end{lemma}
\emph{Proof:} Take $K\in\mathscr{T}_h$, let $\overline{n}$ be the outward normal to $\partial K$, and momentarily assume that $w\in H^3(K)$. An application of integration by parts gives us
\begin{equation}\label{ibp:id}
\langle D^2w,D^2v_h\rangle_K+\langle\Delta w,\nabla v_h\cdot\overline{n}\rangle_{\partial K}-\langle\nabla(\nabla w\cdot\overline{n}),\nabla v_h\rangle_{\partial K}+\langle(\nabla w)^T,\nabla n^T\nabla v_h\rangle_{\partial K} = \langle\Delta w,\Delta v_h\rangle_K.
\end{equation}
A density argument shows that~(\ref{ibp:id}) holds for $w\in H^s(K)$, $s>5/2$.

Applying the identities in~(\ref{identities}) to~(\ref{ibp:id}), and summing over $K\in\mathscr{T}_h$, noting that the normal is constant on faces in $\mathscr{E}^i_h$, we obtain
\begin{equation*}
\begin{split}
&\sum_{K\in\mathscr{T}_h}\langle D^2w,D^2v_h\rangle_K+\sum_{F\in\mathscr{E}^{i,b}_h}\int_F\jump{(\operatorname{div}_\Ta \nabla_\Ta w)(\nabla v_h\cdot n_F)-\nabla_\Ta (\nabla v\cdot n_F)\cdot\nabla_\Ta v_h}\,ds\\
&\,\,\,\,+\sum_{F\in\mathscr{E}^b_h}[\mathcal{H}_F\frac{\partial w}{\partial n_F}\frac{\partial v_h}{\partial n_F}+(\nabla w)^T\nabla n^T_F\nabla v_h]\,ds=\sum_{K\in\mathscr{T}_h}\langle\Delta w,\Delta v_h\rangle_K.
\end{split}
\end{equation*}
The remainder of the argument follows identically as in the proof of Lemma 4 in~\cite{MR3077903}, noting the following calculation, which gives us the final two terms present in~(\ref{Bh*def}):
\begin{equation*}
\begin{split}
(\nabla w)^T\nabla n_F^T\nabla v_h & = ((\nabla_\Ta w)^T+\frac{\partial w}{\partial n_F}n_F^T)(\nabla n_F^T)(\nabla_ Tv_h+\frac{\partial v_h}{\partial n_F}n_F)\\
& = (\nabla_\Ta w)^T\nabla n_F^T\nabla_\Ta v_h+\frac{\partial w}{\partial n_F}(n_F^T\nabla n_F^T\nabla_\Ta v_h)\\
&\,\,\,\,+(\nabla_\Ta w)^T(\nabla n_F^Tn_F)\frac{\partial v_h}{\partial n_F}+\frac{\partial w}{\partial n_F}\frac{\partial v_h}{\partial n_F}n_F^T\nabla n^T_Fn_F\\
& = (\nabla_\Ta w)^T\nabla n_F^T\nabla_\Ta v_h+\frac{\partial w}{\partial n_F}(n_F^T\nabla n_F^T\nabla_\Ta v_h)\\
&\,\,\,\,+(\nabla_\Ta w)^T\left(\frac{1}{2}\nabla |n_F|^2\right)\frac{\partial v_h}{\partial n_F}+\frac{\partial w}{\partial n_F}\frac{\partial v_h}{\partial n_F}n_F^T\left(\frac{1}{2}\nabla|n^T_F|^2\right)
\end{split}
\end{equation*}
Since $|n_F|=1$, it follows that the last two terms above are both zero (and thus, so is their sum).$\quad\quad\square$

The following corollary shows that the method is consistent, that is, if the true solution, $u$, of~(\ref{1}) is sufficiently smooth then $u$ also satisfies~(\ref{method}).
\begin{corollary}
Let $\Omega$ be a piecewise $C^{2,\alpha}$ domain, with $\alpha\in(0,1)$, and let $\Th$ be an exact (possibly curved) simplicial or parallelepipedal mesh on $\Omega$. Assume that $u\in H^2(\Omega)\cap H^1_0(\Omega)$ is a solution of~(\ref{1}).
If $u\in H^s(\Omega;\Th)$, $s>5/2$, then $u$ satisfies
\begin{equation}
A_h(u,v_h)=\sum_{K\in\Th}\langle\gamma f,\Delta v_h\rangle_K\quad\forall w_h\in V_{h,p}.
\end{equation}
\end{corollary}
\emph{Proof:} This follows simply by noting that $u$ satisfies 
\begin{equation*}
\left\{
\begin{aligned}
\gamma Lu&=\gamma f,\quad\mbox{a.e in }\Omega,\\
u &= 0,\quad~~\mbox{on }\partial\Omega,
\end{aligned}
\right.
\end{equation*}
as well as the regularity assumptions necessary for Lemma~\ref{cons:lemma} to hold.$\quad\quad\square$
\subsection{Error estimates}\label{Error estimates}
\begin{theorem}\label{error:est:thm}
Let $\Omega$ be a piecewise $C^{3,\alpha}$-piecewise convex domain, with $\alpha\in(0,1)$. Moreover, assume that $\partial\Omega$ is piecewise $C^m$ for some $m\ge 3$, and let $\{\Th\}_h$ be a regular sequence of simplicial or parallelepipedal meshes with curved faces satisfying~(\ref{meshcond1})-(\ref{Fcontained}) for each $h$. Assume that the sequence of meshes consists of meshes that are regular of order $m$. 
Let $u\in H^2(\Omega)\cap H^1_0(\Omega)$ be the unique solution of~(\ref{1}). Assume that $u\in H^{\mathbf s}(\Omega;\Th)$ with $m>s_K>5/2$ for each $K\in\Th$. Let $c_{\operatorname{stab}}$, $c_{\mathcal{H}}$, $\mu_F$, and $\eta_F$ be chosen as in Theorem~\ref{mainthm} for all $F\in\mathscr{E}^{i,b}_h$, and let $\eta_F>0$ for each $F\in\mathscr{E}^{i,b}_h$. Then, there exists 
a positive constant $C$ independent of $h$ and $u$, but depending on $\max_Ks_K$, such that for the unique solution $u_h$ of~(\ref{method}), we have
\begin{equation}\label{errorest1}
\|u-u_h\|_{h,1}^2\le C\sum_{K\in\Th}h_K^{2t_K-4}\|u\|^2_{H^{s_K}(K)},
\end{equation}
where $t_K=\min(p+1,s_K)$ for each $K\in\Th$.

Note that for the special case of quasi-uniform meshes and uniform polynomial degrees, if $u\in H^s(\Omega)$ with $s>5/2$, the a priori estimate~(\ref{errorest1}) simplifies to 
\begin{equation*}
\|u-u_h\|_{h,1}\le Ch^{\min(p+1,s)-2}\|u\|_{H^s(\Omega)}.
\end{equation*}
Therefore, the convergence rates are optimal with respect to the mesh size.
\end{theorem}
\emph{Proof:} The proof is analogous to the proof of Theorem 9 in~\cite{MR3077903}, Section 5.
It is noteworthy that the proof relies on the existence of a $z_h\in V_{h,p}$ and a constant $C$, independent of $u$, $h_K$ and $p$, but dependent on $\max_Ks_K$, such that for each $K\in\Th$, each nonnegative integer $q\le\min\{s_K,m+1\}$, and each multi-index $\beta$ with $|\beta|<s_K-1/2$, we have
\begin{equation}\label{opt:interp}
\begin{split}
\|u-z_h\|_{H^q(K)}&\le Ch_K^{t_K-q}\|u\|_{H^{s_K}(K)},\\
\|D^\beta(u-z_h)\|_{L^2(\partial K)}&\le Ch^{t_K-|\beta|-1/2}\|u\|_{H^{s_K}(K)}.
\end{split}
\end{equation}
The existence of such a $z_h$ follows from~\cite{babuvska1987hp} and~\cite{MR1014883}. 

The error estimates given by the first inequality in~(\ref{opt:interp}) is given in~\cite{babuvska1987hp} in the context of meshes consisting of simplices and parallelepipeds that do not have curved faces. These results, however, still hold when elements of the mesh are curved. First one must note that the first inequality in~(\ref{opt:interp}) follows from the trace inequality, followed by an application of the second inequality in~(\ref{opt:interp}). Furthermore, in~\cite{MR1014883}, the second bound in~(\ref{opt:interp}) is derived (see Corollary 4.1 in~\cite{MR1014883}).

In order to generalise the proof found in~\cite{MR3077903} to the framework of this paper, it is sufficient to show that for 
$$\xi_h:= z_h-u,$$
and
$$\psi_h:=z_h-u_h$$
we have
\begin{equation}\label{thm:goal}
\begin{split}
&\sum_{F\in\mathscr{E}^b_h}\left\langle\mathcal{H}\frac{\partial\xi_h}{\partial n_F},\frac{\partial \psi_h}{\partial n_F}\right\rangle_F+\langle\nabla_\Ta\xi_h,\nabla n^T_F\nabla_\Ta\psi_h\rangle_F+\left\langle\frac{\partial\xi_h}{\partial n_F},n_F^T\nabla n^T_F\nabla_\Ta\psi_h\right\rangle_F-\left\langle\nabla_\Ta \frac{\partial \psi_h}{\partial n_F},\nabla_\Ta \xi_h\right\rangle_F\\
&~~~~~~~~~~~~~~~~~~~~~~~~~~~\le C\left(\sum_{K\in\Th}h_K^{2t_K-4}\|u\|^2_{H^{s_K}(K)}\right)^{1/2}\|\psi_h\|_{h,1}.
\end{split}
\end{equation}
To establish this bound, we first note that, for any $F\in\Eb$, estimate~(\ref{gradTdn:bound}) also holds for $\psi_h$, that is,
$$\left|\gradTa\left(\frac{\partial \psi_h}{\partial n_F}\right)\right|\le(d-1)(|D^2\psi_h|+|\nabla n_F||\gradTa\psi_h|).$$
From this, we obtain the following:
\begin{equation*}
\begin{split}
&\sum_{F\in\mathscr{E}^b_h}\left\langle\mathcal{H}_F\frac{\partial\xi_h}{\partial n_F},\frac{\partial \psi_h}{\partial n_F}\right\rangle_F+\langle\nabla_\Ta\xi_h,\nabla n^T_F\nabla_\Ta\psi_h\rangle_F+\left\langle\frac{\partial\xi_h}{\partial n_F},n_F^T\nabla n^T_F\nabla_\Ta\psi_h\right\rangle_F-\left\langle\nabla_\Ta \frac{\partial \psi_h}{\partial n_F},\nabla_\Ta \xi_h\right\rangle_F\\
&\le\sum_{F\in\Eb}\left[\left\|\mc^{1/2}_F\frac{\partial\xi_h}{\partial n_F}\right\|_{L^2(F)}\left\|\mc^{1/2}_F\frac{\partial\psi_h}{\partial n_F}\right\|_{L^2(F)}\right.\\
&~~~~~~+\max_{F\in\Eb}\|\nabla n_F^T\|_{L^\infty(F)}(\|\gradTa\xi_h\|_{L^2(F)}\|\gradTa\psi_h\|_{L^2(F)}+\left\|\frac{\partial\xi_h}{\partial n_F}\right\|_{L^2(F)}\|\gradTa\psi_h\|_{L^2(F)})\\
&~~~~~~~~~~~~~+\left.(d-1)\|\gradTa\xi_h\|_{L^2(F)}(\|D^2\psi_h\|_{L^2(F)}+\max_{F\in\Eb}\|\nabla n_F^T\|_{L^\infty(F)}\|\gradTa\psi_h\|_{L^2(F)})\right],
\end{split}
\end{equation*}
and since the quantities
$$\max_{F\in\Eb}\|\nabla n_F^T\|_{L^\infty(F)},\,\max_{F\in\Eb}\|\mc^{1/2}_F\|_{L^\infty(F)}$$
are bounded independently of the mesh size and polynomial degree, we obtain, after an application of the Cauchy--Schwarz inequality for $n$-dimensional vectors:
\begin{equation*}
\begin{split}
&\sum_{F\in\mathscr{E}^b_h}\left\langle\mathcal{H}_F\frac{\partial\xi_h}{\partial n_F},\frac{\partial \psi_h}{\partial n_F}\right\rangle_F+\langle\nabla_\Ta\xi_h,\nabla n^T_F\nabla_\Ta\psi_h\rangle_F+\left\langle\frac{\partial\xi_h}{\partial n_F},n_F^T\nabla n^T_F\nabla_\Ta\psi_h\right\rangle_F-\left\langle\nabla_\Ta \frac{\partial \psi_h}{\partial n_F},\nabla_\Ta \xi_h\right\rangle_F\\
&~~~~~~~~\le C\left(\sum_{F\in\Eb}\left\|\frac{\partial\xi_h}{\partial n_F}\right\|_{L^2(F)}^2+\left(1+\frac{1}{\tilde{h}_F}\right)\|\gradTa\xi_h\|_{L^2(F)}^2\right)^{1/2}\\
&~~~~~~~~~~~~~~~~~~\times\left(\sum_{F\in\Eb}\left\|\mc^{1/2}_F\frac{\partial\psi_h}{\partial n_F}\right\|_{L^2(F)}^2+\|\gradTa\psi_h\|_{L^2(F)}^2+\tilde{h}_F\|D^2\psi_h\|_{L^2(F)}^2\right)^{1/2}.
\end{split}
\end{equation*}
After applying  the inverse inequality~(\ref{bnd:inv:est}) with $m=2$, we obtain
\begin{equation*}
\begin{split}
&\sum_{F\in\Eb}\left\|\mc^{1/2}_F\frac{\partial\psi_h}{\partial n_F}\right\|_{L^2(F)}^2+\|\gradTa\psi_h\|_{L^2(F)}^2+\tilde{h}_F\|D^2\psi_h\|_{L^2(F)}^2\\
&~~~~~~~~~~~\le \sum_{F\in\Eb}\left\|\mc^{1/2}_F\frac{\partial\psi_h}{\partial n_F}\right\|_{L^2(F)}^2+\|\gradTa\psi_h\|_{L^2(F)}^2+C\sum_{K\in\Th}\frac{1}{\tilde{h}_F^{3}}\|\psi_h\|_{L^2(K)}^2\\
&~~~~~~~~~~~\le C\|\psi_h\|_{h,1}^2.
\end{split}
\end{equation*}
We then apply the second interpolation estimate in~(\ref{opt:interp}), yielding,
\begin{equation*}
\begin{split}
\sum_{F\in\Eb}\left\|\frac{\partial\xi_h}{\partial n_F}\right\|_{L^2(F)}^2+(1+\frac{1}{\tilde{h}_F})\|\gradTa\xi_h\|_{L^2(F)}^2&\le C\sum_{F\in\Eb}\frac{1}{\tilde{h}_F}\|\nabla\xi_h\|_{L^2(F)}^2\\
&\le C\sum_{K\in\Th}h_K^{2t_K-4}\|u\|_{H^{s_K}(K)}^2.
\end{split}
\end{equation*}
Combining these two estimates, we obtain~(\ref{thm:goal}).$\quad\quad\square$

\subsection{Quadratic domain approximation} 
In order to prove error estimate~(\ref{errorest1}), and the consistency result~(\ref{cons:statement}), it was required to assume that the triangulations, $\Th$, under consideration approximate the domain exactly. In the case that the domain, $\Omega$, is convex and polytopal, this can be achieved using standard quasi--uniform meshes. 

When the domain is a piecewise $C^{2,\alpha}$ boundary, $\alpha\in(0,1)$, with at least one curved boundary portion, the approach is not so simple. In theory, one can construct exact meshes by considering the boundary, $\partial\Omega$, as a $(d-1)$ dimensional hyper--surface, using the $C^{2,\alpha}$ functions $\{\psi_i\}_{i\in I}$, where $\psi_i:\mathbb R^{d}\to\partial\Omega$ that locally describe the boundary (note that the index set $I$ is determined by $\Omega$).

In practice, this turns out to be somewhat difficult, so instead, one can approximate each map $\psi_i$, $i\in I$, by interpolating it into a Lagrange finite element space, $\mathbb L$, consisting of $d$-dimensional, vector-valued finite element functions.

To define the space $\mathbb L$, we generate the polytopal domain $\Omega_1$ by placing a collection of quasi-uniformly spaced points on $\partial\Omega$, and taking the closed convex hull of these points.  We then take a quasi-uniform triangulation of $\Omega_1$, which we call ${\Th}_{,1}$. Note that in triangulating $\Omega_1$, by ${\Th}_{,1}$, we generate a collection of \emph{affine} maps $\tilde{F}_K:\hat{K}\to\tilde{K}$ of the form~(\ref{def:affine-map}). We then define $\mathbb L$ as follows:
$$\mathbb L:=\{v\in C^0(\overline{\Omega}_1;\mathbb R^d):v|_{\tilde{K}}\in\mathbb P^2(\tilde{K}),\,\forall \tilde{K}\in{\Th}_{,1}\}.$$
We then take the function $x=(x_1,\ldots,x_d)$ which is the coordinate map for the triangulation. Since $x\in\mathbb L$, it admits the representation
$$x = \sum_{j=1}^{N_{\operatorname{int}}}\alpha_j\phi_j+\sum_{j=N_{\operatorname{int}}+1}^{N_{\operatorname{int}}+N_{\partial\Omega_1}}\alpha_j\phi_j,$$
where $\alpha_j\in\Reals$, and the set $\{\phi_j\}_{j=1}^{N_{\operatorname{int}}+N_{\partial\Omega_1}}$ forms a basis of $\mathbb L$, ordered so that $\phi_1,\ldots,\phi_{N_{\operatorname{int}}}$ make up the basis functions associated with the internal degrees of freedom, $x_1,\ldots,x_{N_{\operatorname{int}}}$, and $\phi_{N_{\operatorname{int}}+1},\ldots,$\\$\phi_{N_{\operatorname{int}}+N_{\partial\Omega_1}}$ make up the basis functions associated with the boundary degrees of freedom, $x_{N_{\operatorname{int}}+1},\ldots,$\\$x_{N_{\operatorname{int}}+N_{\partial\Omega_1}}$, that lie on $\partial\Omega_1$.

We then let $$x_h = \sum_{j=1}^{N_{\operatorname{int}}+N_{\partial\Omega_1}}\beta_j\phi_j,$$
with
\begin{equation}\label{betajs}
\beta_j:=\left\{
\begin{aligned}
\alpha_j, & \quad j = 1,\ldots,N_{\operatorname{int}},\\
\psi_i(x_j),& \quad j=N_{\operatorname{int}}+1,\ldots,N_{\operatorname{int}}+N_{\partial\Omega_1},
\end{aligned}
\right.
\end{equation}
for some $i\in I$ (that is, the value of $x_h$ at a given degree of freedom of the finite element space $\mathbb L$ can only differ from $x$ if the degree of freedom lies on $\partial\Omega_1$).

Finally, we define the collection of maps $\Psi_K:\hat{K}\to\mathbb R^d$ by
$$\Psi_K(x_1,\ldots,x_d):=\left\{
\begin{aligned}
\tilde{F}_K & \quad\mbox{if at most one vertex of }\tilde{K}\,\,\mbox{lies on}\,\,\partial\Omega,\\
\pi_h(x_h)& \quad\mbox{otherwise},
\end{aligned}
\right.
$$
where $\pi_h$ is the interpolation operator for the finite element space $\mathbb L$.

From this, we obtain a collection of maps $F_K$ of the form:
$$F_K = \tilde{F}_K+\Phi_K,$$
where $\Phi_K = \Psi_K-\tilde{F}_K$.

Since $\mathbb L$ consists of vector-valued finite element functions of (up to) quadratic order, we obtain a quadratic approximation of $\Omega$. In experiments~\ref{theexp:1},~\ref{theexp:3}, and~\ref{exp:cons} the domain $\Omega$ is the unit disk; in this case we consider the collection $\{\psi_i\}_{i\in I}$, where $I$ is just a singleton set, and the map $\psi_1=\psi$ is given by
$$\psi(x) = x/|x|.$$
In experiment~\ref{theexp:3}, the domain $\Omega$ is the ``key-hole" shaped domain given by~(\ref{Key-hole}); in this case, the collection $\{\psi_i\}_{i\in I}=\{\psi_1,\psi_2\}=\{x/|x|,x\}$, and the choice of $\psi_i$ in~(\ref{betajs}) is determined by whether the degree of freedom lies on a flat or curved boundary portion of $\partial\Omega$.
\begin{remark}
In our experiments, we allow the polynomial degree, $p$, of the finite element space, $V_{h,p}$, vary from $2$ to $4$. For $p>4$, we observe that the quadratic domain approximation becomes dominant, yielding rates of convergence lower than one would expect, were we not committing a so-called ``variational crime" (see~\cite{MR2373954}).
\end{remark}
\end{section}
\begin{section}{Experiments}\label{Experiments}
In this section, we test the robustness of the scheme~(\ref{method}), with the computational domain $\Omega$ taken to be the unit disk, and consider various elliptic operators, $L$, that satisfy the Cordes condition~(\ref{cordes}). In each case, we see that the convergence rates are of the expected order in the various broken Sobolev norms considered, ands in particular in the $\|\cdot\|_{h,1}$--norm, for which we have proven the error bound~(\ref{errorest1}).
\subsection{Experiment 1}\label{theexp:1}
In this experiment, we consider the following problem
\begin{equation}\label{exp:1}
\left\{
\begin{aligned}
\Delta u &= f,\quad\mbox{in}\quad\Omega,\\
u & = 0,\quad\mbox{on}\quad\partial\Omega,\\
\end{aligned}
\right.
\end{equation}
where $\Omega=\{(x,y)\in\mathbb R^2:|x|<1\}$. In this case $$\gamma := \frac{\operatorname{Tr}(A)}{|A|^2}=\frac{\operatorname{Tr}(I_d)}{|I_d|^2}=\frac{I_d:I_d}{I_d:I_d}=1,$$ and the solution of~(\ref{exp:1}) is given by $$u(x,y) = \frac{1}{4}\sin(\pi(x^2+y^2)).$$ In this experiment, we successively increase the degree, $p$, of the finite element space $V_{h,p}$ from $2$ to $4$, and for each fixed degree we refine the mesh quasi--uniformly.
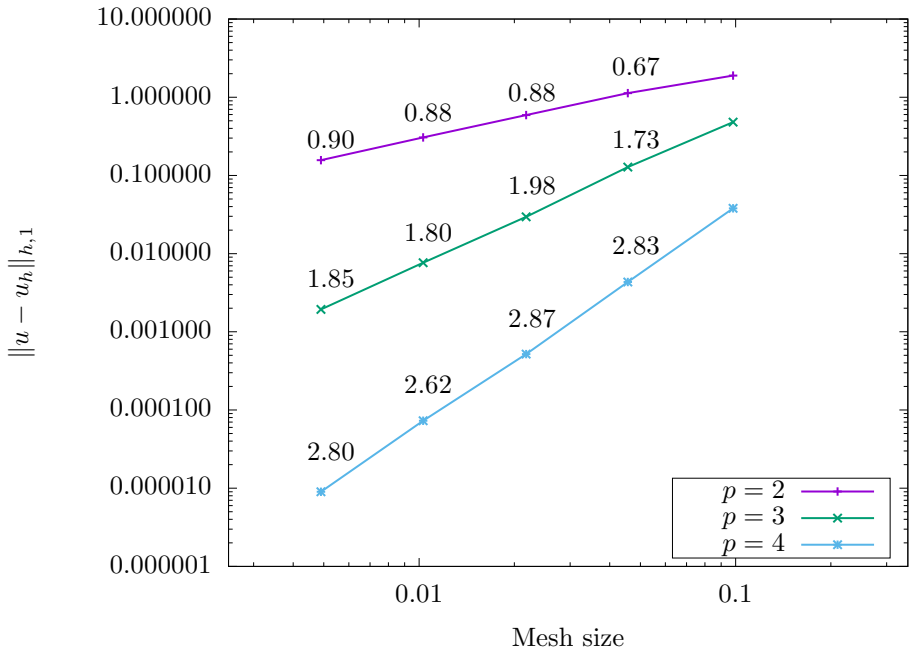
\begin{figure}[h]\label{fig:exp:1}
\input{exp_1_new.tex}
\caption{Convergence rates for the numerical scheme applied to problem~(\ref{exp:1}). The error $\|u-u_h\|_{h,1}$ is plotted against the mesh size $h$ for polynomial degrees ranging from $p=2
$ to $p=4$. We also provide the order of convergence observed.}
\end{figure}
\newpage
\subsection{Experiment 2}\label{theexp:2}
In this experiment, we consider the following problem
\begin{equation}\label{exp:2}
\left\{
\begin{aligned}
\sum_{i,j=1}^d(1+\delta_{ij})\frac{x_i}{|x_i|}\frac{x_j}{|x_j|} D^2_{ij}u&= f,\quad\mbox{in}\quad\Omega,\\
u & = 0,\quad\mbox{on}\quad\partial\Omega,\\
\end{aligned}
\right.
\end{equation}
where $\Omega=\{(x,y)\in\mathbb R^2:|x|<1\}$. In this case $$\gamma = \frac{\operatorname{Tr}(A)}{|A|^2}= \frac{2+x^2/|x|^2+y^2/|y|^2}{8+2x^2y^2/(|x|^2|y|^2)} = 2/5,$$ and $f$ is chosen so that the solution of~(\ref{exp:2}) is given by $$u(x,y) = \frac{1}{4}\sin(\pi(x^2+y^2)).$$ In this experiment, we successively increase the degree, $p$, of the finite element space $V_{h,p}$ from $2$ to $4$, and for each fixed degree we refine the mesh quasi--uniformly.
\begin{figure}[h]\label{fig:exp:2}
\input{exp_2_new.tex}
\caption{Convergence rates for the numerical scheme applied to problem~(\ref{exp:2}). The error $\|u-u_h\|_{H^2(\Omega;\Th)}$ is plotted aganist the mesh size $h$ for polynomial degrees ranging from $p=2
$ to $p=4$. We also provide the order of convergence observed.}
\end{figure}
\newpage
\subsection{Experiment 3}\label{theexp:3}
In this experiment, we consider the PDE given by~(\ref{exp:2}). In this case $f$ is chosen so that the solution of~(\ref{exp:2}) is given by 
$$u(x,y) = \frac{1}{4}\sin(\pi(x^2+y^2)).$$
We have also taken the $\Omega$ to be ``key-hole" shaped domain
$$\{x^2+y^2<1:y\ge1/\sqrt{2}\}\cup[-1/\sqrt{2},1/\sqrt{2}]\times[-3,1/\sqrt{2}],$$
thus demonstrating the applicability of our numerical method on piecewise curved nonconvex domains. Furthermore, the boundary value problem considered is inhomogeneous. In order to extend our numerical method~(\ref{method}) to this case, we simply modify the right hand side as follows (denoting $g$ to be the restriction of $u$ the boundary, $\partial\Omega$)
\begin{equation*}
\begin{aligned}
A_h(u_h,v_h) & = \sum_{K\in\Th}\langle\gamma f,\Delta v_h\rangle_K+\sum_{F\in\Eb}[\mu_f\langle\nabla_\Ta g,\nabla_\Ta v_h\rangle_F+\eta_F\langle g,v_h\rangle_F]\\
&-\frac{1}{2}\sum_{F\in\Eb}\langle\operatorname{div}_\Ta\nabla_\Ta g,\nabla v_h\cdot n_F\rangle_F+\langle\nabla_\Ta(\nabla v_h\cdot n_F),\nabla_\Ta g\rangle_F].
\end{aligned}
\end{equation*}
In this experiment, we successively increase the degree, $p$, of the finite element space $V_{h,p}$ from $2$ to $4$, and for each fixed degree we refine the mesh quasi--uniformly.
\begin{figure}[h]\label{fig:exp:3}
\input{exp_3_new.tex}
\caption{Convergence rates for the numerical scheme applied to problem~(\ref{exp:2}), on a nonconvex domain given by~(\ref{Key-hole}). The normalised error values $\|u-u_h\|_{h,1}/\|u\|_{h,1}$ are plotted against the mesh size $h$ with the polynomial degree $p=2$. for polynomial degrees ranging from $p=2
$ to $p=4$. We also provide the order of convergence observed.}
\end{figure}
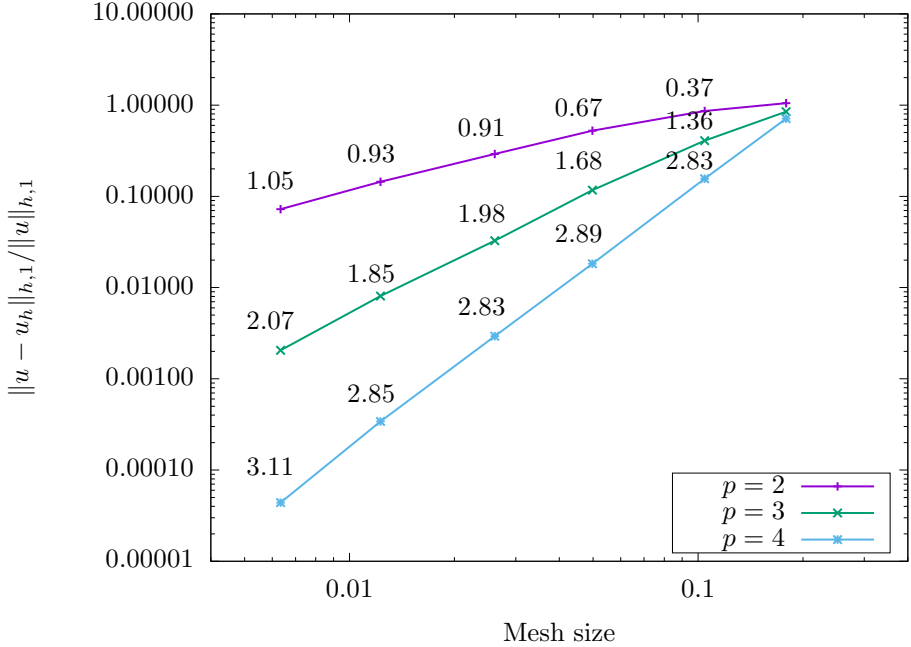
\begin{figure}
\begin{center}
\scalebox{0.6}{
\begin{tikzpicture}
 \draw[ultra thick] (1+0.5,-2) arc (-45:225:2cm);
 \draw[ultra thick] (1+0.5,-2) -- (1+0.5,-2-3-1-0.7071067811865);
 \draw[ultra thick] (1+0.5-2.8284271247461,-2) -- (1+0.5-2.8284271247461,-2-3-1-0.7071067811865);
 \draw[ultra thick] (1+0.5-2.8284271247461,-2-3-1-0.7071067811865) -- (1+0.5,-2-3-1-0.7071067811865);
 \draw[ultra thick] (6+0.5+2,-2-1.1767766952966+1) arc (90:270:1cm);
 \draw[ultra thick] (6+0.5+2,-1-1-1.1767766952966+1) -- (6+0.5+2,-1-1.1767766952966+1);
 \draw[ultra thick] (6+0.5+2,-1-1.1767766952966+1) -- (2+0.5+2,-1-1.1767766952966+1);
 \draw[ultra thick] (2+0.5+2,-1-1.1767766952966+1) --  (2+0.5+2,-1-1.1767766952966-4+1);
 \draw[ultra thick] (2+0.5+2,-1-1.1767766952966-4+1) -- (6+0.5+2,-1-1.1767766952966-4+1);
 \draw[ultra thick] (6+0.5+2,-1-1.1767766952966-4+1) -- (6+0.5+2,-1-1.1767766952966-4+1+1);
     \end{tikzpicture}
    }
\end{center}
\caption{Examples of the ``key-hole" shaped domain (left) given by~(\ref{Key-hole}), and a domain with a boundary portion of strictly negative curvature (right) given by~(\ref{diskremoved}).}\label{key:hole:mesh}
\end{figure}
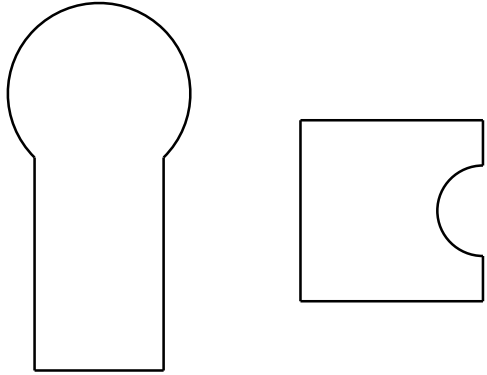

\newpage
\subsection{Experiment 4 - Consistency}\label{exp:cons}
As mentioned in the introduction, the bilinear form $B_{h,*}$ defined by~(\ref{Bh*def}) includes terms that are necessary for the consistency of the method, arising from the curvature of the boundary. These terms are not present in the method presented in~\cite{MR3077903}, and the following experiment shows the necessity of including these new terms; in particular, we see both a lack of consistency, and error results inferior to those produced by the new method~(\ref{method}). 

In the results that follow, we provide the consistency residual
$$\operatorname{Res}(w):=B_{h,*}(w,w)-\sum_{K\in\Th}\langle\Delta w,\Delta w\rangle_K,$$
for a function $w\in H^s(\Omega;\Th)\cap H^2(\Omega)\cap H^1_0(\Omega)\cap V_{h,p},$ $s>5/2$, which numerically validates Lemma~\ref{cons:lemma}, 
as well as the error results arising from one mesh refinement. In the first set of results, we implement the method presented in~\cite{MR3077903}, which we shall call ``Method A", for problem~(\ref{exp:2}), with a quadratic approximation space. In the second set, we implement the method presented in this paper, which we shall call ``Method B", for the same problem, with the same approximation space.
%

\begin{table}[H]
\begin{center}
\begin{tabular}{|c|c|c|c|c|c|c|}\hline
\multicolumn{1}{|c|}{Method} &
\multicolumn{1}{|c|}{Refinement no.} &
\multicolumn{1}{|c|}{ $\operatorname{Res}(w)$} &
\multicolumn{2}{|c|}{$\|u-u_h\|_{L^2(\Omega)}$, EOC} & 
\multicolumn{2}{|c|}{$\|u-u_h\|_{h,1}$, EOC} \\ \hline
A & 1 & -6.244 & 1.534e-01 & - & 9.597e-01 & - \\ \hline
A & 2 & -6.273 & 1.468e-01 & 5.735e-02 & 9.286e-01 & 4.304e-02\\ \hline
B & 1 & -2.692e-05 & 2.538e-04 & - & 8.810e-02 & - \\ \hline
B & 2 & -1.754e-06 & 3.7801e-05 & 2.48 & 3.176e-02 & 1.33  \\ \hline
\end{tabular}
\end{center}
\end{table}
\begin{remark}[Polynomial domain approximation]
We have assumed that the mesh of the computational domain is exact, i.e., that~(\ref{Exact}) holds. In practice, we are able to preserve optimal (in the sense of~(\ref{errorest1})) error bounds, by using a polynomial approximation of the domain; when the polynomial degree of the domain approximation matches that of the finite element space, it is referred to as ``isoparametric approximation" (see~\cite{MR1014883}).
\end{remark}
\end{section}
\begin{section}{Conclusion}\label{Conclusion}
We have extended the framework introduced in~\cite{MR3077903}, allowing for domains with curved boundaries. We have tested the robustness of this new method (given by~(\ref{method})) with numerical experiments involving elliptic operators with discontinuous coefficients, on a uniformly convex domain that has a curved boundary, and a nonconvex domain with both flat and curved boundary portions with strictly positive curvature. Furthermore, experiment~\ref{exp:cons} validated the necessity of the modifications to the method found in~\cite{MR3077903}, that are present in our new method~(\ref{method}).

For the two computational domains considered, in order to verify the error estimates present in Section~\ref{Analysis of the numerical method} we used meshes consisting of curved triangles with edges were defined by a combination of polynomial and affine mappings. It would be an interesting avenue for future research to consider ellipsoidal and oval--shaped domains, and domains with a boundaries that are not piecewise $C^\infty$, and to see what happens in cases of largely varying curvature.

The type of problems under consideration (problems in nondivergence form on curved domains) pose many analytical and computational difficulties, whilst housing a large variety of applications; in this paper we have developed a method that produces optimal error results. This inference has been validated by the analysis in Section~\ref{Analysis of the numerical method}, and the numerical experiments found in Section~\ref{Experiments}.
\end{section}
\bibliographystyle{plain}
\bibliography{toskaweckiBIB}

\end{document}

%% file: exp_1_new.tex
\begingroup
  \makeatletter
  \providecommand\color[2][]{%
    \GenericError{(gnuplot) \space\space\space\@spaces}{%
      Package color not loaded in conjunction with
      terminal option `colourtext'%
    }{See the gnuplot documentation for explanation.%
    }{Either use 'blacktext' in gnuplot or load the package
      color.sty in LaTeX.}%
    \renewcommand\color[2][]{}%
  }%
  \providecommand\includegraphics[2][]{%
    \GenericError{(gnuplot) \space\space\space\@spaces}{%
      Package graphicx or graphics not loaded%
    }{See the gnuplot documentation for explanation.%
    }{The gnuplot epslatex terminal needs graphicx.sty or graphics.sty.}%
    \renewcommand\includegraphics[2][]{}%
  }%
  \providecommand\rotatebox[2]{#2}%
  \@ifundefined{ifGPcolor}{%
    \newif\ifGPcolor
    \GPcolortrue
  }{}%
  \@ifundefined{ifGPblacktext}{%
    \newif\ifGPblacktext
    \GPblacktexttrue
  }{}%
  \let\gplgaddtomacro\g@addto@macro
  \gdef\gplbacktext{}%
  \gdef\gplfronttext{}%
  \makeatother
  \ifGPblacktext
    \def\colorrgb#1{}%
    \def\colorgray#1{}%
  \else
    \ifGPcolor
      \def\colorrgb#1{\color[rgb]{#1}}%
      \def\colorgray#1{\color[gray]{#1}}%
      \expandafter\def\csname LTw\endcsname{\color{white}}%
      \expandafter\def\csname LTb\endcsname{\color{black}}%
      \expandafter\def\csname LTa\endcsname{\color{black}}%
      \expandafter\def\csname LT0\endcsname{\color[rgb]{1,0,0}}%
      \expandafter\def\csname LT1\endcsname{\color[rgb]{0,1,0}}%
      \expandafter\def\csname LT2\endcsname{\color[rgb]{0,0,1}}%
      \expandafter\def\csname LT3\endcsname{\color[rgb]{1,0,1}}%
      \expandafter\def\csname LT4\endcsname{\color[rgb]{0,1,1}}%
      \expandafter\def\csname LT5\endcsname{\color[rgb]{1,1,0}}%
      \expandafter\def\csname LT6\endcsname{\color[rgb]{0,0,0}}%
      \expandafter\def\csname LT7\endcsname{\color[rgb]{1,0.3,0}}%
      \expandafter\def\csname LT8\endcsname{\color[rgb]{0.5,0.5,0.5}}%
    \else
      \def\colorrgb#1{\color{black}}%
      \def\colorgray#1{\color[gray]{#1}}%
      \expandafter\def\csname LTw\endcsname{\color{white}}%
      \expandafter\def\csname LTb\endcsname{\color{black}}%
      \expandafter\def\csname LTa\endcsname{\color{black}}%
      \expandafter\def\csname LT0\endcsname{\color{black}}%
      \expandafter\def\csname LT1\endcsname{\color{black}}%
      \expandafter\def\csname LT2\endcsname{\color{black}}%
      \expandafter\def\csname LT3\endcsname{\color{black}}%
      \expandafter\def\csname LT4\endcsname{\color{black}}%
      \expandafter\def\csname LT5\endcsname{\color{black}}%
      \expandafter\def\csname LT6\endcsname{\color{black}}%
      \expandafter\def\csname LT7\endcsname{\color{black}}%
      \expandafter\def\csname LT8\endcsname{\color{black}}%
    \fi
  \fi
    \setlength{\unitlength}{0.0500bp}%
    \ifx\gptboxheight\undefined%
      \newlength{\gptboxheight}%
      \newlength{\gptboxwidth}%
      \newsavebox{\gptboxtext}%
    \fi%
    \setlength{\fboxrule}{0.5pt}%
    \setlength{\fboxsep}{1pt}%
\begin{picture}(7200.00,5040.00)%
    \gplgaddtomacro\gplbacktext{%
      \csname LTb\endcsname
      \put(1606,704){\makebox(0,0)[r]{\strut{}0.000001}}%
      \put(1606,1292){\makebox(0,0)[r]{\strut{}0.000010}}%
      \put(1606,1880){\makebox(0,0)[r]{\strut{}0.000100}}%
      \put(1606,2468){\makebox(0,0)[r]{\strut{}0.001000}}%
      \put(1606,3055){\makebox(0,0)[r]{\strut{}0.010000}}%
      \put(1606,3643){\makebox(0,0)[r]{\strut{}0.100000}}%
      \put(1606,4231){\makebox(0,0)[r]{\strut{}1.000000}}%
      \put(1606,4819){\makebox(0,0)[r]{\strut{}10.000000}}%
      \put(3159,484){\makebox(0,0){\strut{}$0.01$}}%
      \put(5519,484){\makebox(0,0){\strut{}$0.1$}}%
    }%
    \gplgaddtomacro\gplfronttext{%
      \csname LTb\endcsname
      \put(198,2761){\rotatebox{-270}{\makebox(0,0){\strut{}$\|u-u_h\|_{h,1}$}}}%
      \put(4270,154){\makebox(0,0){\strut{}Mesh size}}%
      \put(4775,4450){\makebox(0,0){\strut{}$0.67$}}
      \put(4775,3900){\makebox(0,0){\strut{}$1.73$}}
      \put(4775,3100){\makebox(0,0){\strut{}$2.83$}}%
      \put(4000,4250){\makebox(0,0){\strut{}$0.88$}}
      \put(4000,3550){\makebox(0,0){\strut{}$1.98$}}
      \put(4000,2550){\makebox(0,0){\strut{}$2.87$}}%
      \put(3225,4100){\makebox(0,0){\strut{}$0.88$}}
      \put(3225,3200){\makebox(0,0){\strut{}$1.80$}}
      \put(3225,2050){\makebox(0,0){\strut{}$2.62$}}%
      \put(2500,3900){\makebox(0,0){\strut{}$0.90$}}
      \put(2500,2850){\makebox(0,0){\strut{}$1.85$}}
      \put(2500,1550){\makebox(0,0){\strut{}$2.80$}}%
      \csname LTb\endcsname
      \put(5888,1267){\makebox(0,0)[r]{\strut{}$~~~~p=2$}}%
      \csname LTb\endcsname
      \put(5888,1067){\makebox(0,0)[r]{\strut{}$~~~~p=3$}}%
      \csname LTb\endcsname
      \put(5888,867){\makebox(0,0)[r]{\strut{}$~~~~p=4$}}%
    }%
    \gplbacktext
    \put(0,0){\includegraphics{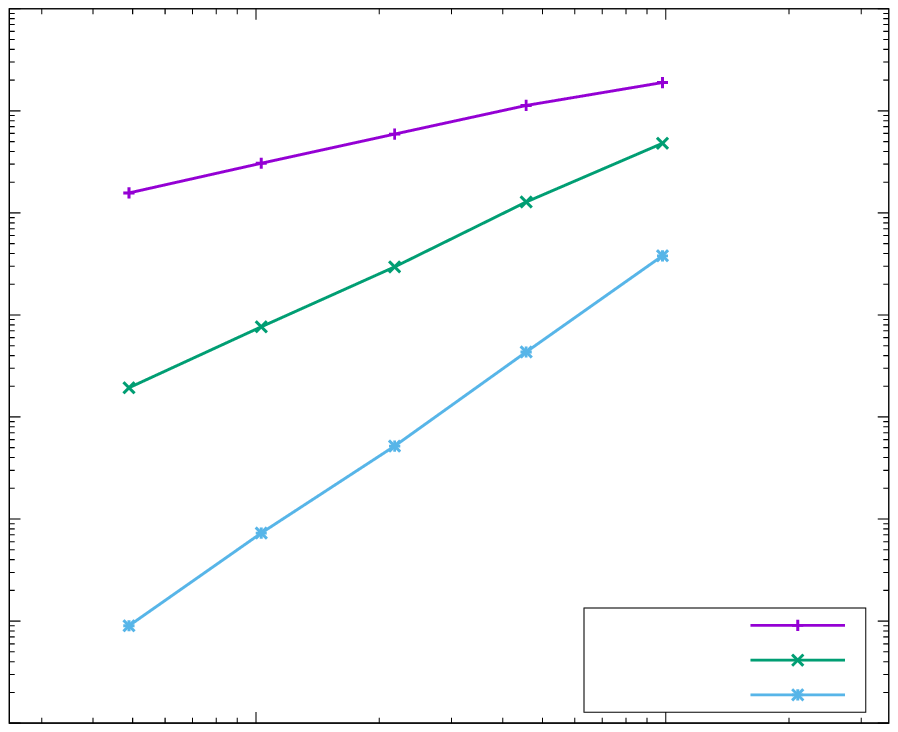}}%
    \gplfronttext
  \end{picture}%
\endgroup

%% file: exp_2_new.tex
\begingroup
  \makeatletter
  \providecommand\color[2][]{%
    \GenericError{(gnuplot) \space\space\space\@spaces}{%
      Package color not loaded in conjunction with
      terminal option `colourtext'%
    }{See the gnuplot documentation for explanation.%
    }{Either use 'blacktext' in gnuplot or load the package
      color.sty in LaTeX.}%
    \renewcommand\color[2][]{}%
  }%
  \providecommand\includegraphics[2][]{%
    \GenericError{(gnuplot) \space\space\space\@spaces}{%
      Package graphicx or graphics not loaded%
    }{See the gnuplot documentation for explanation.%
    }{The gnuplot epslatex terminal needs graphicx.sty or graphics.sty.}%
    \renewcommand\includegraphics[2][]{}%
  }%
  \providecommand\rotatebox[2]{#2}%
  \@ifundefined{ifGPcolor}{%
    \newif\ifGPcolor
    \GPcolortrue
  }{}%
  \@ifundefined{ifGPblacktext}{%
    \newif\ifGPblacktext
    \GPblacktexttrue
  }{}%
  \let\gplgaddtomacro\g@addto@macro
  \gdef\gplbacktext{}%
  \gdef\gplfronttext{}%
  \makeatother
  \ifGPblacktext
    \def\colorrgb#1{}%
    \def\colorgray#1{}%
  \else
    \ifGPcolor
      \def\colorrgb#1{\color[rgb]{#1}}%
      \def\colorgray#1{\color[gray]{#1}}%
      \expandafter\def\csname LTw\endcsname{\color{white}}%
      \expandafter\def\csname LTb\endcsname{\color{black}}%
      \expandafter\def\csname LTa\endcsname{\color{black}}%
      \expandafter\def\csname LT0\endcsname{\color[rgb]{1,0,0}}%
      \expandafter\def\csname LT1\endcsname{\color[rgb]{0,1,0}}%
      \expandafter\def\csname LT2\endcsname{\color[rgb]{0,0,1}}%
      \expandafter\def\csname LT3\endcsname{\color[rgb]{1,0,1}}%
      \expandafter\def\csname LT4\endcsname{\color[rgb]{0,1,1}}%
      \expandafter\def\csname LT5\endcsname{\color[rgb]{1,1,0}}%
      \expandafter\def\csname LT6\endcsname{\color[rgb]{0,0,0}}%
      \expandafter\def\csname LT7\endcsname{\color[rgb]{1,0.3,0}}%
      \expandafter\def\csname LT8\endcsname{\color[rgb]{0.5,0.5,0.5}}%
    \else
      \def\colorrgb#1{\color{black}}%
      \def\colorgray#1{\color[gray]{#1}}%
      \expandafter\def\csname LTw\endcsname{\color{white}}%
      \expandafter\def\csname LTb\endcsname{\color{black}}%
      \expandafter\def\csname LTa\endcsname{\color{black}}%
      \expandafter\def\csname LT0\endcsname{\color{black}}%
      \expandafter\def\csname LT1\endcsname{\color{black}}%
      \expandafter\def\csname LT2\endcsname{\color{black}}%
      \expandafter\def\csname LT3\endcsname{\color{black}}%
      \expandafter\def\csname LT4\endcsname{\color{black}}%
      \expandafter\def\csname LT5\endcsname{\color{black}}%
      \expandafter\def\csname LT6\endcsname{\color{black}}%
      \expandafter\def\csname LT7\endcsname{\color{black}}%
      \expandafter\def\csname LT8\endcsname{\color{black}}%
    \fi
  \fi
    \setlength{\unitlength}{0.0500bp}%
    \ifx\gptboxheight\undefined%
      \newlength{\gptboxheight}%
      \newlength{\gptboxwidth}%
      \newsavebox{\gptboxtext}%
    \fi%
    \setlength{\fboxrule}{0.5pt}%
    \setlength{\fboxsep}{1pt}%
\begin{picture}(7200.00,5040.00)%
    \gplgaddtomacro\gplbacktext{%
      \csname LTb\endcsname
      \put(1606,704){\makebox(0,0)[r]{\strut{}0.000001}}%
      \put(1606,1292){\makebox(0,0)[r]{\strut{}0.000010}}%
      \put(1606,1880){\makebox(0,0)[r]{\strut{}0.000100}}%
      \put(1606,2468){\makebox(0,0)[r]{\strut{}0.001000}}%
      \put(1606,3055){\makebox(0,0)[r]{\strut{}0.010000}}%
      \put(1606,3643){\makebox(0,0)[r]{\strut{}0.100000}}%
      \put(1606,4231){\makebox(0,0)[r]{\strut{}1.000000}}%
      \put(1606,4819){\makebox(0,0)[r]{\strut{}10.000000}}%
      \put(3159,484){\makebox(0,0){\strut{}$0.01$}}%
      \put(5519,484){\makebox(0,0){\strut{}$0.1$}}%
    }%
    \gplgaddtomacro\gplfronttext{%
      \csname LTb\endcsname
      \put(198,2761){\rotatebox{-270}{\makebox(0,0){\strut{}$\|u-u_h\|_{h,1}$}}}%
      \put(4270,154){\makebox(0,0){\strut{}Mesh size}}%
      \put(4775,4450){\makebox(0,0){\strut{}$0.68$}}
      \put(4775,3900){\makebox(0,0){\strut{}$1.73$}}
      \put(4775,3100){\makebox(0,0){\strut{}$2.83$}}%
      \put(4000,4250){\makebox(0,0){\strut{}$0.88$}}
      \put(4000,3550){\makebox(0,0){\strut{}$1.98$}}
      \put(4000,2550){\makebox(0,0){\strut{}$2.87$}}%
      \put(3225,4100){\makebox(0,0){\strut{}$0.88$}}
      \put(3225,3200){\makebox(0,0){\strut{}$1.80$}}
      \put(3225,2050){\makebox(0,0){\strut{}$2.62$}}%
      \put(2500,3900){\makebox(0,0){\strut{}$0.90$}}
      \put(2500,2850){\makebox(0,0){\strut{}$1.85$}}
      \put(2500,1550){\makebox(0,0){\strut{}$2.80$}}%

      \csname LTb\endcsname
      \put(5888,1267){\makebox(0,0)[r]{\strut{}$~~~~p=2$}}%
      \csname LTb\endcsname
      \put(5888,1067){\makebox(0,0)[r]{\strut{}$~~~~p=3$}}%
      \csname LTb\endcsname
      \put(5888,867){\makebox(0,0)[r]{\strut{}$~~~~p=4$}}%
    }%
    \gplbacktext
    \put(0,0){\includegraphics{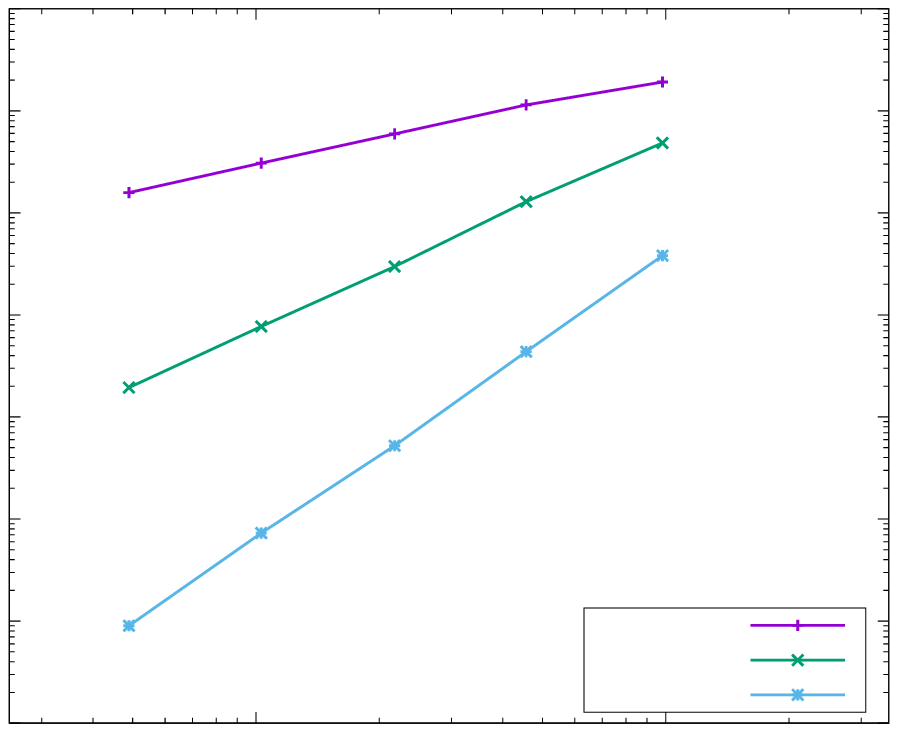}}%
    \gplfronttext
  \end{picture}%
\endgroup

%% file: exp_3_new.tex
\begingroup
  \makeatletter
  \providecommand\color[2][]{%
    \GenericError{(gnuplot) \space\space\space\@spaces}{%
      Package color not loaded in conjunction with
      terminal option `colourtext'%
    }{See the gnuplot documentation for explanation.%
    }{Either use 'blacktext' in gnuplot or load the package
      color.sty in LaTeX.}%
    \renewcommand\color[2][]{}%
  }%
  \providecommand\includegraphics[2][]{%
    \GenericError{(gnuplot) \space\space\space\@spaces}{%
      Package graphicx or graphics not loaded%
    }{See the gnuplot documentation for explanation.%
    }{The gnuplot epslatex terminal needs graphicx.sty or graphics.sty.}%
    \renewcommand\includegraphics[2][]{}%
  }%
  \providecommand\rotatebox[2]{#2}%
  \@ifundefined{ifGPcolor}{%
    \newif\ifGPcolor
    \GPcolortrue
  }{}%
  \@ifundefined{ifGPblacktext}{%
    \newif\ifGPblacktext
    \GPblacktexttrue
  }{}%
  \let\gplgaddtomacro\g@addto@macro
  \gdef\gplbacktext{}%
  \gdef\gplfronttext{}%
  \makeatother
  \ifGPblacktext
    \def\colorrgb#1{}%
    \def\colorgray#1{}%
  \else
    \ifGPcolor
      \def\colorrgb#1{\color[rgb]{#1}}%
      \def\colorgray#1{\color[gray]{#1}}%
      \expandafter\def\csname LTw\endcsname{\color{white}}%
      \expandafter\def\csname LTb\endcsname{\color{black}}%
      \expandafter\def\csname LTa\endcsname{\color{black}}%
      \expandafter\def\csname LT0\endcsname{\color[rgb]{1,0,0}}%
      \expandafter\def\csname LT1\endcsname{\color[rgb]{0,1,0}}%
      \expandafter\def\csname LT2\endcsname{\color[rgb]{0,0,1}}%
      \expandafter\def\csname LT3\endcsname{\color[rgb]{1,0,1}}%
      \expandafter\def\csname LT4\endcsname{\color[rgb]{0,1,1}}%
      \expandafter\def\csname LT5\endcsname{\color[rgb]{1,1,0}}%
      \expandafter\def\csname LT6\endcsname{\color[rgb]{0,0,0}}%
      \expandafter\def\csname LT7\endcsname{\color[rgb]{1,0.3,0}}%
      \expandafter\def\csname LT8\endcsname{\color[rgb]{0.5,0.5,0.5}}%
    \else
      \def\colorrgb#1{\color{black}}%
      \def\colorgray#1{\color[gray]{#1}}%
      \expandafter\def\csname LTw\endcsname{\color{white}}%
      \expandafter\def\csname LTb\endcsname{\color{black}}%
      \expandafter\def\csname LTa\endcsname{\color{black}}%
      \expandafter\def\csname LT0\endcsname{\color{black}}%
      \expandafter\def\csname LT1\endcsname{\color{black}}%
      \expandafter\def\csname LT2\endcsname{\color{black}}%
      \expandafter\def\csname LT3\endcsname{\color{black}}%
      \expandafter\def\csname LT4\endcsname{\color{black}}%
      \expandafter\def\csname LT5\endcsname{\color{black}}%
      \expandafter\def\csname LT6\endcsname{\color{black}}%
      \expandafter\def\csname LT7\endcsname{\color{black}}%
      \expandafter\def\csname LT8\endcsname{\color{black}}%
    \fi
  \fi
    \setlength{\unitlength}{0.0500bp}%
    \ifx\gptboxheight\undefined%
      \newlength{\gptboxheight}%
      \newlength{\gptboxwidth}%
      \newsavebox{\gptboxtext}%
    \fi%
    \setlength{\fboxrule}{0.5pt}%
    \setlength{\fboxsep}{1pt}%
\begin{picture}(7200.00,5040.00)%
    \gplgaddtomacro\gplbacktext{%
      \csname LTb\endcsname
      \put(1474,704){\makebox(0,0)[r]{\strut{}0.00001}}%
      \put(1474,1390){\makebox(0,0)[r]{\strut{}0.00010}}%
      \put(1474,2076){\makebox(0,0)[r]{\strut{}0.00100}}%
      \put(1474,2762){\makebox(0,0)[r]{\strut{}0.01000}}%
      \put(1474,3447){\makebox(0,0)[r]{\strut{}0.10000}}%
      \put(1474,4133){\makebox(0,0)[r]{\strut{}1.00000}}%
      \put(1474,4819){\makebox(0,0)[r]{\strut{}10.00000}}%
      \put(2640,484){\makebox(0,0){\strut{}$0.01$}}%
      \put(5239,484){\makebox(0,0){\strut{}$0.1$}}%
    }%
    \gplgaddtomacro\gplfronttext{%
      \csname LTb\endcsname
      \put(198,2761){\rotatebox{-270}{\makebox(0,0){\strut{}$\|u-u_h\|_{h,1}/\|u\|_{h,1}$}}}%
      \put(4204,154){\makebox(0,0){\strut{}Mesh size}}%
      \put(5175,4250){\makebox(0,0){\strut{}$0.37$}}
      \put(5175,4000){\makebox(0,0){\strut{}$1.36$}}
      \put(5175,3700){\makebox(0,0){\strut{}$2.83$}}%
      \put(4350,4100){\makebox(0,0){\strut{}$0.67$}}
      \put(4350,3700){\makebox(0,0){\strut{}$1.68$}}
      \put(4350,3150){\makebox(0,0){\strut{}$2.89$}}%
      \put(3625,3950){\makebox(0,0){\strut{}$0.91$}}
      \put(3625,3300){\makebox(0,0){\strut{}$1.98$}}
      \put(3625,2600){\makebox(0,0){\strut{}$2.83$}}%
      \put(2800,3750){\makebox(0,0){\strut{}$0.93$}}
      \put(2800,2850){\makebox(0,0){\strut{}$1.85$}}
      \put(2800,1950){\makebox(0,0){\strut{}$2.85$}}%
      \put(2050,3550){\makebox(0,0){\strut{}$1.05$}}
      \put(2050,2500){\makebox(0,0){\strut{}$2.07$}}
      \put(2050,1400){\makebox(0,0){\strut{}$3.11$}}%
      \csname LTb\endcsname
      \put(5888,1267){\makebox(0,0)[r]{\strut{}$~~~~p=2$}}%
      \csname LTb\endcsname
      \put(5888,1067){\makebox(0,0)[r]{\strut{}$~~~~p=3$}}%
      \csname LTb\endcsname
      \put(5888,867){\makebox(0,0)[r]{\strut{}$~~~~p=4$}}%
    }%
    \gplbacktext
    \put(0,0){\includegraphics{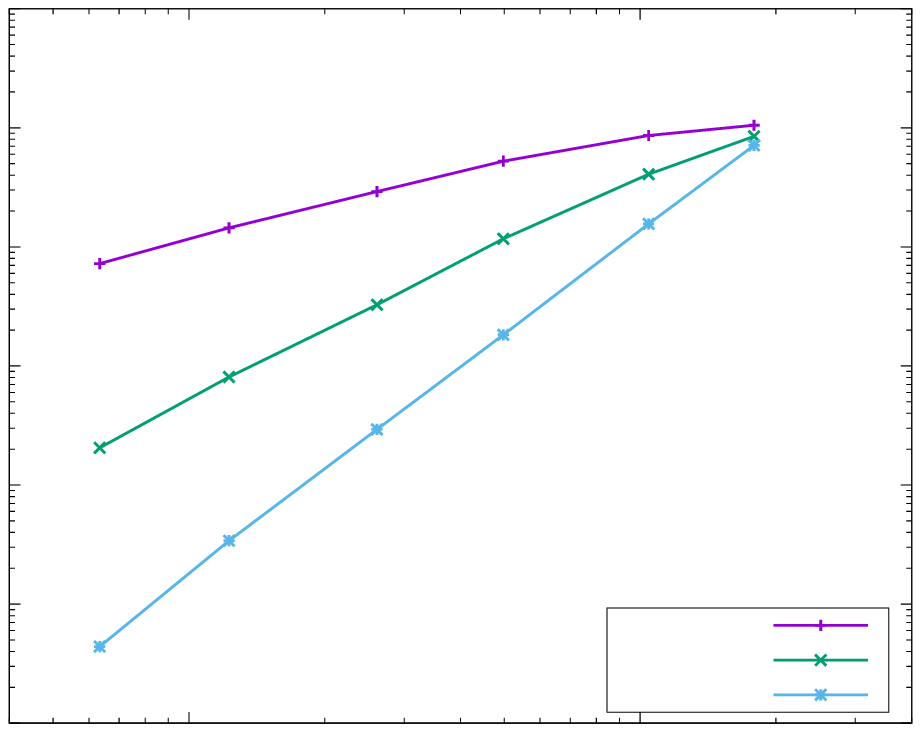}}%
        \gplfronttext
  \end{picture}%
\endgroup

%% file: DGFEM-Linear-Dirichlet-Curved.bbl
\begin{thebibliography}{10}

\bibitem{babuvska1987hp}
I.~Babu{\v{s}}ka and M.~Suri.
\newblock The $ hp $ version of the finite element method with quasiuniform
  meshes.
\newblock {\em RAIRO-Mod{\'e}lisation math{\'e}matique et analyse
  num{\'e}rique}, 21(2):199--238, 1987.

\bibitem{MR1014883}
C.~Bernardi.
\newblock Optimal finite-element interpolation on curved domains.
\newblock {\em SIAM J. Numer. Anal.}, 26(5):1212--1240, 1989.

\bibitem{MR2916369}
S.~C. Brenner and M.~Neilan.
\newblock Finite element approximations of the three dimensional
  {M}onge-{A}mp\`ere equation.
\newblock {\em ESAIM Math. Model. Numer. Anal.}, 46(5):979--1001, 2012.

\bibitem{MR2373954}
S.~C. Brenner and L.~R. Scott.
\newblock {\em The mathematical theory of finite element methods}, volume~15 of
  {\em Texts in Applied Mathematics}.
\newblock Springer, New York, third edition, 2008.

\bibitem{MR1426885}
L.~A. Caffarelli.
\newblock Boundary regularity of maps with convex potentials. {II}.
\newblock {\em Ann. of Math. (2)}, 144(3):453--496, 1996.

\bibitem{MR2179357}
W.~H. Fleming and H.~M. Soner.
\newblock {\em Controlled {M}arkov processes and viscosity solutions},
  volume~25 of {\em Stochastic Modelling and Applied Probability}.
\newblock Springer, New York, second edition, 2006.

\bibitem{MR1814364}
D.~Gilbarg and N.~S. Trudinger.
\newblock {\em Elliptic partial differential equations of second order}.
\newblock Classics in Mathematics. Springer-Verlag, Berlin, 2001.
\newblock Reprint of the 1998 edition.

\bibitem{MR3033005}
M.~Jensen and I.~Smears.
\newblock On the convergence of finite element methods for
  {H}amilton-{J}acobi-{B}ellman equations.
\newblock {\em SIAM J. Numer. Anal.}, 51(1):137--162, 2013.

\bibitem{MR2260015}
A.~Maugeri, D.~K. Palagachev, and L.~G. Softova.
\newblock {\em Elliptic and parabolic equations with discontinuous
  coefficients}, volume 109 of {\em Mathematical Research}.
\newblock Wiley-VCH Verlag Berlin GmbH, Berlin, 2000.

\bibitem{MR3162358}
M.~Neilan.
\newblock Finite element methods for fully nonlinear second order {PDE}s based
  on a discrete {H}essian with applications to the {M}onge-{A}mp\`ere equation.
\newblock {\em J. Comput. Appl. Math.}, 263:351--369, 2014.

\bibitem{MR3653852}
M.~Neilan, A.~J. Salgado, and W.~Zhang.
\newblock Numerical analysis of strongly nonlinear {PDE}s.
\newblock {\em Acta Numer.}, 26:137--303, 2017.

\bibitem{MR0180763}
A.~V. Pogorelov.
\newblock {\em Monge-{A}mp\`ere equations of elliptic type}.
\newblock Translated from the first Russian edition by Leo F. Boron with the
  assistance of Albert L. Rabenstein and Richard C. Bollinger. P. Noordhoff,
  Ltd., Groningen, 1964.

\bibitem{MR3077903}
I.~Smears and E.~S{\"u}li.
\newblock Discontinuous {G}alerkin finite element approximation of
  nondivergence form elliptic equations with {C}ord\`es coefficients.
\newblock {\em SIAM J. Numer. Anal.}, 51(4):2088--2106, 2013.

\bibitem{MR3196952}
I.~Smears and E.~S{\"u}li.
\newblock Discontinuous {G}alerkin finite element approximation of
  {H}amilton-{J}acobi-{B}ellman equations with {C}ord\`es coefficients.
\newblock {\em SIAM J. Numer. Anal.}, 52(2):993--1016, 2014.

\bibitem{MR1454261}
J.~Urbas.
\newblock On the second boundary value problem for equations of
  {M}onge-{A}mp\`ere type.
\newblock {\em J. Reine Angew. Math.}, 487:115--124, 1997.

\bibitem{MR1624426}
J.~Urbas.
\newblock Oblique boundary value problems for equations of {M}onge-{A}mp\`ere
  type.
\newblock {\em Calc. Var. Partial Differential Equations}, 7(1):19--39, 1998.

\end{thebibliography}
